%
%

\def\Date{2008/05/09}


\ifx\pdfoutput\jamaisdefined\else
\input supp-pdf.tex \pdfoutput=1 \pdfcompresslevel=9

\fi

%

\magnification=1200
\hsize=11.25cm
\vsize=18cm
\parskip 0pt
\parindent=12pt
\voffset=1cm
\hoffset=1cm



\catcode'32=9

\font\tenpc=cmcsc10
\font\eightpc=cmcsc8
\font\eightrm=cmr8
\font\eighti=cmmi8
\font\eightsy=cmsy8
\font\eightbf=cmbx8
\font\eighttt=cmtt8
\font\eightit=cmti8
\font\eightsl=cmsl8
\font\sixrm=cmr6
\font\sixi=cmmi6
\font\sixsy=cmsy6
\font\sixbf=cmbx6

\skewchar\eighti='177 \skewchar\sixi='177
\skewchar\eightsy='60 \skewchar\sixsy='60

\catcode`@=11

\def\tenpoint{%
  \textfont0=\tenrm \scriptfont0=\sevenrm \scriptscriptfont0=\fiverm
  \def\rm{\fam\z@\tenrm}%
  \textfont1=\teni \scriptfont1=\seveni \scriptscriptfont1=\fivei
  \def\oldstyle{\fam\@ne\teni}%
  \textfont2=\tensy \scriptfont2=\sevensy \scriptscriptfont2=\fivesy
  \textfont\itfam=\tenit
  \def\it{\fam\itfam\tenit}%
  \textfont\slfam=\tensl
  \def\sl{\fam\slfam\tensl}%
  \textfont\bffam=\tenbf \scriptfont\bffam=\sevenbf
  \scriptscriptfont\bffam=\fivebf
  \def\bf{\fam\bffam\tenbf}%
  \textfont\ttfam=\tentt
  \def\tt{\fam\ttfam\tentt}%
  \abovedisplayskip=12pt plus 3pt minus 9pt
  \abovedisplayshortskip=0pt plus 3pt
  \belowdisplayskip=12pt plus 3pt minus 9pt
  \belowdisplayshortskip=7pt plus 3pt minus 4pt
  \smallskipamount=3pt plus 1pt minus 1pt
  \medskipamount=6pt plus 2pt minus 2pt
  \bigskipamount=12pt plus 4pt minus 4pt
  \normalbaselineskip=12pt
  \setbox\strutbox=\hbox{\vrule height8.5pt depth3.5pt width0pt}%
  \let\bigf@ntpc=\tenrm \let\smallf@ntpc=\sevenrm
  \let\petcap=\tenpc
  \normalbaselines\rm}

\def\eightpoint{%
  \textfont0=\eightrm \scriptfont0=\sixrm \scriptscriptfont0=\fiverm
  \def\rm{\fam\z@\eightrm}%
  \textfont1=\eighti \scriptfont1=\sixi \scriptscriptfont1=\fivei
  \def\oldstyle{\fam\@ne\eighti}%
  \textfont2=\eightsy \scriptfont2=\sixsy \scriptscriptfont2=\fivesy
  \textfont\itfam=\eightit
  \def\it{\fam\itfam\eightit}%
  \textfont\slfam=\eightsl
  \def\sl{\fam\slfam\eightsl}%
  \textfont\bffam=\eightbf \scriptfont\bffam=\sixbf
  \scriptscriptfont\bffam=\fivebf
  \def\bf{\fam\bffam\eightbf}%
  \textfont\ttfam=\eighttt
  \def\tt{\fam\ttfam\eighttt}%
  \abovedisplayskip=9pt plus 2pt minus 6pt
  \abovedisplayshortskip=0pt plus 2pt
  \belowdisplayskip=9pt plus 2pt minus 6pt
  \belowdisplayshortskip=5pt plus 2pt minus 3pt
  \smallskipamount=2pt plus 1pt minus 1pt
  \medskipamount=4pt plus 2pt minus 1pt
  \bigskipamount=9pt plus 3pt minus 3pt
  \normalbaselineskip=9pt
  \setbox\strutbox=\hbox{\vrule height7pt depth2pt width0pt}%
  \let\bigf@ntpc=\eightrm \let\smallf@ntpc=\sixrm
  \let\petcap=\eightpc
  \normalbaselines\rm}
\catcode`@=12

\tenpoint


\long\def\irmaaddress{{%
\bigskip
\eightpoint
\rightline{\quad
\vtop{\halign{\hfil##\hfil\cr
I.R.M.A. UMR 7501\cr
Universit\'e Louis Pasteur et CNRS,\cr
7, rue Ren\'e-Descartes\cr
F-67084 Strasbourg, France\cr
{\tt guoniu@math.u-strasbg.fr}\cr}}\quad}
}}


\font\tengoth=eufm10
\def\goth#1{\hbox{\tengoth #1}}


\catcode`\@=11
\def\pc#1#2|{{\bigf@ntpc #1\penalty \@MM\hskip\z@skip\smallf@ntpc%
	\uppercase{#2}}}
\catcode`\@=12

\def\pointir{\discretionary{.}{}{.\kern.35em---\kern.7em}\nobreak
   \hskip 0em plus .3em minus .4em }

\def\qed{\quad\raise -2pt\hbox{\vrule\vbox to 10pt{\hrule width 4pt
   \vfill\hrule}\vrule}}

\def\rem#1|{\par\medskip{{\it #1}\pointir}}

\def\vspace[#1]{\noalign{\vskip#1}}

\def\abstract#1{\vbox{\eightpoint\narrower\narrower 
\pc ABSTRACT|\pointir #1}}


\def\section#1{\goodbreak\par\vskip .3cm\centerline{\bf #1}
   \par\nobreak\vskip 3pt }

\long\def\th#1|#2\endth{\par\medbreak
   {\petcap #1\pointir}{\it #2}\par\medbreak}

\def\article#1|#2|#3|#4|#5|#6|#7|
    {{\leftskip=7mm\noindent
     \hangindent=2mm\hangafter=1
     \llap{{\tt [#1]}\hskip.35em}{\petcap#2}\pointir
     #3, {\sl #4}, {\bf #5} ({\oldstyle #6}),
     pp.\nobreak\ #7.\par}}
\def\livre#1|#2|#3|#4|
    {{\leftskip=7mm\noindent
    \hangindent=2mm\hangafter=1
    \llap{{\tt [#1]}\hskip.35em}{\petcap#2}\pointir
    {\sl #3}, #4.\par}}
\def\divers#1|#2|#3|
    {{\leftskip=7mm\noindent
    \hangindent=2mm\hangafter=1
     \llap{{\tt [#1]}\hskip.35em}{\petcap#2}\pointir
     #3.\par}}



\catcode`\@=11
\def\c@rr@#1{\vbox{%
  \hrule height \ep@isseur%
   \hbox{\vrule width\ep@isseur\vbox to \t@ille{%
           \vfil\hbox  to \t@ille{\hfil#1\hfil}\vfil}%
            \vrule width\ep@isseur}%
      \hrule height \ep@isseur}}
\def\ytableau#1#2#3#4{\vbox{%
  \gdef\ep@isseur{#2}
   \gdef\t@ille{#1}
    \def\\##1{\c@rr@{$#3 ##1$}}
  \lineskiplimit=-30cm \baselineskip=\t@ille%
    \advance \baselineskip by \ep@isseur%
     \halign{%
      \hfil$##$\hfil&&\kern -\ep@isseur%
       \hfil$##$\hfil \crcr#4\crcr}}}%
\catcode`\@=12

\def\Grille{\setbox13=\vbox to 5mm{\hrule width 110mm\vfill}
\setbox13=\vbox{\offinterlineskip
   \copy13\copy13\copy13\copy13\copy13\copy13\copy13\copy13
   \copy13\copy13\copy13\copy13\box13\hrule width 110mm}
\setbox14=\hbox to 5mm{\vrule height 65mm\hfill}
\setbox14=\hbox{\copy14\copy14\copy14\copy14\copy14\copy14
   \copy14\copy14\copy14\copy14\copy14\copy14\copy14\copy14
   \copy14\copy14\copy14\copy14\copy14\copy14\copy14\copy14\box14}
\ht14=0pt\dp14=0pt\wd14=0pt
\setbox13=\vbox to 0pt{\vss\box13\offinterlineskip\box14}
\wd13=0pt\box13}


\def\fleche(#1,#2)\dir(#3,#4)\long#5{%
\noalign{\nointerlineskip\leftput(#1,#2){\vector(#3,#4){#5}}\nointerlineskip}}


\def\hfl#1#2#3{\smash{\mathop{\hbox to#3{\rightarrowfill}}\limits
^{\scriptstyle#1}_{\scriptstyle#2}}}

\def\gfl#1#2#3{\smash{\mathop{\hbox to#3{\leftarrowfill}}\limits
^{\scriptstyle#1}_{\scriptstyle#2}}}


 \message{`lline' & `vector' macros from LaTeX}
 \catcode`@=11
\def\{{\relax\ifmmode\lbrace\else$\lbrace$\fi}
\def\}{\relax\ifmmode\rbrace\else$\rbrace$\fi}
\def\newcount{\alloc@0\count\countdef\insc@unt}
\def\newdimen{\alloc@1\dimen\dimendef\insc@unt}
\def\newwrite{\alloc@7\write\chardef\sixt@@n}

\newwrite\@unused
\newcount\@tempcnta
\newcount\@tempcntb
\newdimen\@tempdima
\newdimen\@tempdimb
\newbox\@tempboxa

\def\@spaces{\space\space\space\space}
\def\@whilenoop#1{}
\def\@whiledim#1\do #2{\ifdim #1\relax#2\@iwhiledim{#1\relax#2}\fi}
\def\@iwhiledim#1{\ifdim #1\let\@nextwhile=\@iwhiledim
        \else\let\@nextwhile=\@whilenoop\fi\@nextwhile{#1}}
\def\@badlinearg{\@latexerr{Bad \string\line\space or \string\vector
   \space argument}}
\def\@latexerr#1#2{\begingroup
\edef\@tempc{#2}\expandafter\errhelp\expandafter{\@tempc}%
\def\@eha{Your command was ignored.
^^JType \space I <command> <return> \space to replace it
  with another command,^^Jor \space <return> \space to continue without it.}
\def\@ehb{You've lost some text. \space \@ehc}
\def\@ehc{Try typing \space <return>
  \space to proceed.^^JIf that doesn't work, type \space X <return> \space to
  quit.}
\def\@ehd{You're in trouble here.  \space\@ehc}

\typeout{LaTeX error. \space See LaTeX manual for explanation.^^J
 \space\@spaces\@spaces\@spaces Type \space H <return> \space for
 immediate help.}\errmessage{#1}\endgroup}
\def\typeout#1{{\let\protect\string\immediate\write\@unused{#1}}}

\font\tenln    = line10
\font\tenlnw   = linew10

\newdimen\@wholewidth
\newdimen\@halfwidth
\newdimen\unitlength 

\unitlength =1pt


\def\thinlines{\let\@linefnt\tenln \let\@circlefnt\tencirc
  \@wholewidth\fontdimen8\tenln \@halfwidth .5\@wholewidth}
\def\thicklines{\let\@linefnt\tenlnw \let\@circlefnt\tencircw
  \@wholewidth\fontdimen8\tenlnw \@halfwidth .5\@wholewidth}

\def\linethickness#1{\@wholewidth #1\relax \@halfwidth .5\@wholewidth}

\newif\if@negarg

\def\lline(#1,#2)#3{\@xarg #1\relax \@yarg #2\relax
\@linelen=#3\unitlength
\ifnum\@xarg =0 \@vline
  \else \ifnum\@yarg =0 \@hline \else \@sline\fi
\fi}

\def\@sline{\ifnum\@xarg< 0 \@negargtrue \@xarg -\@xarg \@yyarg -\@yarg
  \else \@negargfalse \@yyarg \@yarg \fi
\ifnum \@yyarg >0 \@tempcnta\@yyarg \else \@tempcnta -\@yyarg \fi
\ifnum\@tempcnta>6 \@badlinearg\@tempcnta0 \fi
\setbox\@linechar\hbox{\@linefnt\@getlinechar(\@xarg,\@yyarg)}%
\ifnum \@yarg >0 \let\@upordown\raise \@clnht\z@
   \else\let\@upordown\lower \@clnht \ht\@linechar\fi
\@clnwd=\wd\@linechar
\if@negarg \hskip -\wd\@linechar \def\@tempa{\hskip -2\wd\@linechar}\else
     \let\@tempa\relax \fi
\@whiledim \@clnwd <\@linelen \do
  {\@upordown\@clnht\copy\@linechar
   \@tempa
   \advance\@clnht \ht\@linechar
   \advance\@clnwd \wd\@linechar}%
\advance\@clnht -\ht\@linechar
\advance\@clnwd -\wd\@linechar
\@tempdima\@linelen\advance\@tempdima -\@clnwd
\@tempdimb\@tempdima\advance\@tempdimb -\wd\@linechar
\if@negarg \hskip -\@tempdimb \else \hskip \@tempdimb \fi
\multiply\@tempdima \@m
\@tempcnta \@tempdima \@tempdima \wd\@linechar \divide\@tempcnta \@tempdima
\@tempdima \ht\@linechar \multiply\@tempdima \@tempcnta
\divide\@tempdima \@m
\advance\@clnht \@tempdima
\ifdim \@linelen <\wd\@linechar
   \hskip \wd\@linechar
  \else\@upordown\@clnht\copy\@linechar\fi}

\def\@hline{\ifnum \@xarg <0 \hskip -\@linelen \fi
\vrule height \@halfwidth depth \@halfwidth width \@linelen
\ifnum \@xarg <0 \hskip -\@linelen \fi}

\def\@getlinechar(#1,#2){\@tempcnta#1\relax\multiply\@tempcnta 8
\advance\@tempcnta -9 \ifnum #2>0 \advance\@tempcnta #2\relax\else
\advance\@tempcnta -#2\relax\advance\@tempcnta 64 \fi
\char\@tempcnta}

\def\vector(#1,#2)#3{\@xarg #1\relax \@yarg #2\relax
\@linelen=#3\unitlength
\ifnum\@xarg =0 \@vvector
  \else \ifnum\@yarg =0 \@hvector \else \@svector\fi
\fi}

\def\@hvector{\@hline\hbox to 0pt{\@linefnt
\ifnum \@xarg <0 \@getlarrow(1,0)\hss\else
    \hss\@getrarrow(1,0)\fi}}

\def\@vvector{\ifnum \@yarg <0 \@downvector \else \@upvector \fi}

\def\@svector{\@sline
\@tempcnta\@yarg \ifnum\@tempcnta <0 \@tempcnta=-\@tempcnta\fi
\ifnum\@tempcnta <5
  \hskip -\wd\@linechar
  \@upordown\@clnht \hbox{\@linefnt  \if@negarg
  \@getlarrow(\@xarg,\@yyarg) \else \@getrarrow(\@xarg,\@yyarg) \fi}%
\else\@badlinearg\fi}

\def\@getlarrow(#1,#2){\ifnum #2 =\z@ \@tempcnta='33\else
\@tempcnta=#1\relax\multiply\@tempcnta \sixt@@n \advance\@tempcnta
-9 \@tempcntb=#2\relax\multiply\@tempcntb \tw@
\ifnum \@tempcntb >0 \advance\@tempcnta \@tempcntb\relax
\else\advance\@tempcnta -\@tempcntb\advance\@tempcnta 64
\fi\fi\char\@tempcnta}

\def\@getrarrow(#1,#2){\@tempcntb=#2\relax
\ifnum\@tempcntb < 0 \@tempcntb=-\@tempcntb\relax\fi
\ifcase \@tempcntb\relax \@tempcnta='55 \or
\ifnum #1<3 \@tempcnta=#1\relax\multiply\@tempcnta
24 \advance\@tempcnta -6 \else \ifnum #1=3 \@tempcnta=49
\else\@tempcnta=58 \fi\fi\or
\ifnum #1<3 \@tempcnta=#1\relax\multiply\@tempcnta
24 \advance\@tempcnta -3 \else \@tempcnta=51\fi\or
\@tempcnta=#1\relax\multiply\@tempcnta
\sixt@@n \advance\@tempcnta -\tw@ \else
\@tempcnta=#1\relax\multiply\@tempcnta
\sixt@@n \advance\@tempcnta 7 \fi\ifnum #2<0 \advance\@tempcnta 64 \fi
\char\@tempcnta}

\def\@vline{\ifnum \@yarg <0 \@downline \else \@upline\fi}

\def\@upline{\hbox to \z@{\hskip -\@halfwidth \vrule
  width \@wholewidth height \@linelen depth \z@\hss}}

\def\@downline{\hbox to \z@{\hskip -\@halfwidth \vrule
  width \@wholewidth height \z@ depth \@linelen \hss}}

\def\@upvector{\@upline\setbox\@tempboxa\hbox{\@linefnt\char'66}\raise
     \@linelen \hbox to\z@{\lower \ht\@tempboxa\box\@tempboxa\hss}}

\def\@downvector{\@downline\lower \@linelen
      \hbox to \z@{\@linefnt\char'77\hss}}

\thinlines

\newcount\@xarg
\newcount\@yarg
\newcount\@yyarg
\newcount\@multicnt
\newdimen\@xdim
\newdimen\@ydim
\newbox\@linechar
\newdimen\@linelen
\newdimen\@clnwd
\newdimen\@clnht
\newdimen\@dashdim
\newbox\@dashbox
\newcount\@dashcnt
 \catcode`@=12


\newbox\tbox
\newbox\tboxa

\def\leftzer#1{\setbox\tbox=\hbox to 0pt{#1\hss}%
     \ht\tbox=0pt \dp\tbox=0pt \box\tbox}

\def\rightzer#1{\setbox\tbox=\hbox to 0pt{\hss #1}%
     \ht\tbox=0pt \dp\tbox=0pt \box\tbox}

\def\centerzer#1{\setbox\tbox=\hbox to 0pt{\hss #1\hss}%
     \ht\tbox=0pt \dp\tbox=0pt \box\tbox}

%
\def\image(#1,#2)#3{\vbox to #1{\offinterlineskip
    \vss #3 \vskip #2}}


\def\leftput(#1,#2)#3{\setbox\tboxa=\hbox{%
    \kern #1\unitlength
    \raise #2\unitlength\hbox{\leftzer{#3}}}%
    \ht\tboxa=0pt \wd\tboxa=0pt \dp\tboxa=0pt\box\tboxa}

\def\rightput(#1,#2)#3{\setbox\tboxa=\hbox{%
    \kern #1\unitlength
    \raise #2\unitlength\hbox{\rightzer{#3}}}%
    \ht\tboxa=0pt \wd\tboxa=0pt \dp\tboxa=0pt\box\tboxa}

\def\centerput(#1,#2)#3{\setbox\tboxa=\hbox{%
    \kern #1\unitlength
    \raise #2\unitlength\hbox{\centerzer{#3}}}%
    \ht\tboxa=0pt \wd\tboxa=0pt \dp\tboxa=0pt\box\tboxa}

\unitlength=1mm

\def\cput(#1,#2)#3{\noalign{\nointerlineskip\centerput(#1,#2){#3}
                             \nointerlineskip}}


\ifx\pdfoutput\jamaisdefined
\input epsf

\fi


\parskip 0pt plus 1pt

\def\article#1|#2|#3|#4|#5|#6|#7|
    {{\leftskip=7mm\noindent
     \hangindent=2mm\hangafter=1
     \llap{{\tt [#1]}\hskip.35em}{#2},\quad %
     #3, {\sl #4}, {\bf #5} ({\oldstyle #6}),
     pp.\nobreak\ #7.\par}}
\def\livre#1|#2|#3|#4|
    {{\leftskip=7mm\noindent
    \hangindent=2mm\hangafter=1
    \llap{{\tt [#1]}\hskip.35em}{#2},\quad %
    {\sl #3}, #4.\par}}
\def\divers#1|#2|#3|
    {{\leftskip=7mm\noindent
    \hangindent=2mm\hangafter=1
     \llap{{\tt [#1]}\hskip.35em}{#2},\quad %
     #3.\par}}

\def\l{\lambda}

\def\mod{\mathop{\rm mod}}

\def\lig{\mathop{\goth g}}
\def\tmax{{\max_t}}

\def\enslettre#1{\font\zzzz=msbm10 \hbox{\zzzz #1}}
\def\setZ{\mathop{\enslettre Z}}
\def\setN{\mathop{\enslettre N}}
\def\setP{\mathop{\cal P}}
\def\setS{\mathop{\goth S}}


\rightline{\Date}
\bigskip

\centerline{\bf An explicit expansion formula 
for the powers of the Euler Product}
\centerline{\bf in terms of partition hook lengths}
\bigskip
\centerline{Guo-Niu HAN}
\bigskip\medskip

\abstract{
We discover an explicit expansion formula for the powers $s$ of 
the Euler Product (or Dedekind $\eta$-function) in terms of 
hook lengths of partitions, where the exponent $s$ is {\it any} complex 
number.  Several classical formulas have been derived
for certain integers $s$
by Euler, Jacobi, Klein, Fricke, Atkin, Winquist, Dyson and Macdonald.
In particular, Macdonald obtained expansion formulas 
for the integer exponents $s$ for which there exists a semi-simple Lie algebra 
of dimension~$s$. 
For the type $A_l^{(a)}$ he has expressed the $(t^2-1)$-st power
of the Euler Product as a sum of weighted integer vectors of length $t$
for any integer $t$.  
Kostant has considered the general case for any positive integer $s$ 
and obtained further properties.
\smallskip
The present paper proposes a new approach. We convert 
the weighted vectors of {\it length} $t$ used by Macdonald in his 
identity for type $A_l^{(a)}$ to weighted partitions with 
{\it free parameter} $t$, so that a new identity on the latter 
combinatorial structures can be derived without any restrictions on~$t$.
The surprise is that the 
weighted partitions have a very simple form in terms of hook lengths
of partitions.  
As applications of our formula, we find some 
new identities about hook lengths, including the ``marked hook 
formula". We also improve a result due to Kostant.
The proof of the Main Theorem is based on Macdonald's identity for $A_l^{(a)}$
and on the properties of
a bijection 
between $t$-cores and integer vectors constructed by Garvan, Kim and Stanton.
}

\bigskip

\def\itemm#1|#2| {{\leftskip=7mm\noindent \hangindent=2mm\hangafter=1 
    \llap{\S#1.\ }{#2}\par}}

\def\itemmm#1|#2| {{\leftskip=14mm\noindent \hangindent=2mm\hangafter=1 
    \llap{#1.\ }{#2}\par}}

\centerline{\bf Summary}
\footnote{}{%
\eightpoint
Andrei Okounkov pointed out that the main identity 
already appeared in his joint paper arXiv:hep-th/0306238.
The present paper will remain on arXiv verbatim and not be published 
anywhere else. 
Parts of the results, as well as several new ones
 are reproduced in a forthcoming paper
	arXiv:0805.1398v1 [math.CO].

}

\smallskip

\itemm1|Introduction. The Main Theorem. Selected results.|

\itemm2|Basic consequences and specializations.|

\itemmm2.1|Equivalent forms.|

\itemmm2.2|Corollaries.|

\itemmm2.3|Specialization for $\beta=0$. 
Generating function for partitions.|

\itemmm2.4|Specialization for $\beta=1$.|

\itemmm2.5|Specialization for $\beta=\infty$.
Classical hook length formula and 
the Robinson-Schensted-Knuth correspondence.|

\itemmm2.6|Specialization for $\beta=-1$.|

\itemmm2.7|Specialization for $\beta=2$. 
Euler's pentagonal theorem.
Example for illustrating the Main Theorem.|

\itemmm2.8|Specialization for $\beta=25$. Ramanujan $\tau$-function.
Example for illustrating the Main Theorem.|

\itemm3|Specialization for $\beta=4$. Jacobi's triple product formula.|

\itemm4|Specialization for $\beta=9$.|

\itemm5|Proof of the Main Theorem.|

\itemmm5.1|Fundamental properties of $t$-cores and $V$-codings.|

\itemmm5.2|The bijection $\phi_V$ and an example.|

\itemmm5.3|Proof of the first property.|

\itemmm5.4|Proof of the second property.|

\itemmm5.5|End of the proof of the Main Theorem.|

\itemm6|New formulas about hook lengths.|

\itemmm6.1|Comparing the coefficients of $\beta$.|

\itemmm6.2|Stanley-Elder-Bessenrodt-Bacher-Manivel Theorem.|

\itemmm6.3|Comparing the coefficients of $\beta^2$.|

\itemmm6.4|Comparing the coefficients of $\beta^{n}x^n$ and $\beta^{n-1}x^n$. 
The marked hook formula.|

\itemmm6.5|Comparing the coefficients of $\beta^{n-2}x^n$ and of 
$\beta^{n-3}x^n$.|

\itemm7|Improvement of a result due to Kostant.|
\itemm8|The magic partition formula.|
\itemm9|Reversion of the Euler Product.|


\def\sec{1}
\section{\sec. Introduction} 
The powers of the Euler Product and the hook lengths of partitions 
are two mathematical objects widely studied in the Theory of
Partitions, in Algebraic Combinatorics and Group Representation Theory.
In the present paper we establish a new connection by giving an explicit 
expansion formula for all the powers $s$ of the Euler
Product in terms of partition hook lengths, where the exponent $s$
is any complex number. Recall that the {\it Euler Product}
is the infinite product $\prod_{m\geq 0} (1-x^m)$. A variation of the
Euler Product, called the {\it Dedekind $\eta$-function}, is defined by 
$\eta(x)=x^{1/24} \prod_{m\geq 0} (1-x^m).$
The following two formulas [Eu83; An76, p.11, p.21]
go back to Euler (the pentagonal theorem)
$$
\prod_{m\geq 1}  (1-x^m)
= \sum_{k=-\infty}^{\infty}(-1)^k x^{k(3k+1)/2}
\leqno{(\sec.1)}
$$
and Jacobi (triple product identity)
$$
\prod_{m\geq 1} (1-x^m)^3=\sum_{m\geq 0} (-1)^m (2m+1)x^{m(m+1)/2}.
\leqno{(\sec.2)}
$$
Further explicit formulas for the
powers of the Euler Product
$$
\prod_{m\geq 1} (1-x^m)^s  =  \sum_{k\geq 0} f_k(s) x^k \leqno{(\sec.3})
$$
have been derived for certain integers
$$s=1,3,8,10,14,15,21,24,26,28,35,36,\ldots \leqno{(\sec.4)}$$
by 
Klein and Fricke for $s=8$, 
Atkin for $s=14, 26$,  Winquist for $s=10$, 
and Dyson for $s=24, \ldots$ [Wi69; Dy72].  
The paper entitled ``Affine root systems and Dedekind's $\eta$-function",
written by Macdonald in 1972, is a milestone in the study of powers
of Euler Product [Ma72].  The review of this paper for MathSciNet,
written by Verma [Ve], contains seven pages!  
It has also inspired several followers, see
[Ka74; Mo75; Ko76; Le78; Ko04; Mi85; AF02; CFP05; RS06].
The main achievement of
Macdonald was to unify all the 
well-known formulas for the integers $s$ listed in (\sec.4), 
except for $s=26$. 
He obtained an expansion formula of 
$$\prod_{m\geq 0} (1-x^m)^{\dim \lig}\leqno{(\sec.5)}$$
for every semi-simple Lie algebra $\lig$. 
In the case of type $A_l^{(a)}$, i.e., type $A_l$ with $l$ even [Ma72, p.134],
he expressed the $(t^2-1)$-st power
of the Euler Product as a sum of weighted integer vectors of length $t$
for each odd positive integer $t$:
$$
\eta(x)^{t^2-1} = c_0 \sum_{(v_0,\ldots, v_{t-1})} \prod_{i<j} (v_i-v_j) 
x^{(v_0^2+v_1^2+\cdots+v_{t-1}^2)/(2t)}.
$$
Following this direction it seems difficult 
to obtain more expansion formulas  for other exponents, because $t$ is 
a vector length and has to be an integer.
Kostant considered the general case for positive integer $s$ 
and obtained further properties [Ko04].
\medskip

The present paper proposes a new approach. The main difficulty
is to find an appropriate ``other object" and convert
the weighted vector of {\it length} $t$ in 
Macdonald's identity to a
weighted ``other object" with 
{\it free parameter} $t$.
In fact, we find out
that the ``other object" is merely the classical partition of integer
and that the
weighted partition has a surprisingly simple form in terms of 
hook lengths.
\medskip

Let us describe our Main Theorem.  
The basic notions needed here can be found in 
[Ma95, p.1; St99, p.287;  La01, p.1; Kn98, p.59; An76, p.1].
A {\it partition}~$\l$ is a sequence of positive 
integers $\l=(\l_1, \l_2,\cdots, \l_\ell)$ such that 
$\l_1\geq \l_2 \geq \cdots \geq \l_\ell>0$.
The integers
$(\l_i)_{i=1,2,\ldots, \ell}$ are called the {\it parts} of~$\l$,
the number $\ell$ of parts being the
{\it length} of $\l$ denoted by $\ell(\l)$.  
The sum of its parts $\l_1+ \l_2+\cdots+ \l_\ell$ is
denoted by $|\l|$.
Let $n$ be an integer, a partition 
$\l$ is said to be a partition of $n$ if $|\l|=n$. We write $\l\vdash n$.
The set of all partitions of $n$
is denoted by $\setP(n)$. 
The set of all partitions is denoted by~$\setP$,
so that $$\setP=\bigcup_{n\geq 0} \setP(n).$$
Each partition can be represented by its Ferrers diagram. For example,
$\l=(6,3,3,2)$ is a partition  and its Ferrers diagram is reproduced in 
Fig.~\sec.1.

{
\long\def\maplebegin#1\mapleend{}

\maplebegin

# --------------- begin maple ----------------------

# Copy the following text  to "makefig.mpl"
# then in maple > read("makefig.mpl");
# it will create a file "z_fig_by_maple.tex"

#\unitlength=1pt

Hu:= 12.4; # height quantities
Lu:= Hu; # large unity

X0:=-95.0; Y0:=15.6; # origin position

File:=fopen("z_fig_by_maple.tex", WRITE);

mhook:=proc(x,y,lenx, leny)
local i, d,sp, yy, xx, ct;
	sp:=Hu/8;
	ct:=0;
	for xx from x*Lu+X0 to x*Lu+X0+Hu by sp do
		yy := y*Hu+Y0; 
		fprintf(File, "\\vline(
				xx,   yy+sp*ct-0.2,    Lu*leny-sp*ct+0.1);
		ct:=ct+1;
	od:
	
	ct:=0;
	for yy from y*Hu+Y0 to y*Hu+Y0+Hu by sp do
		xx := x*Lu+X0; 
		fprintf(File, "\\hline(
				xx+sp*ct-0.2,   yy,    Lu*lenx-sp*ct+0.1);
		ct:=ct+1;
	od:

end;

mhook(1,1,2,3);

fclose(File);
# -------------------- end maple -------------------------
\mapleend

\setbox1=\hbox{$
\def\b{\\{\hbox{}}}
\ytableau{12pt}{0.4pt}{}
{\b &\b       \cr 
 \b &\b &\b   \cr
 \b &\b &\b   \cr
 \b &\b &\b &\b &\b &\b  \cr
\noalign{\vskip 3pt}
\noalign{\hbox{Fig.~\sec.1. Partition}}
}$
}
\setbox2=\hbox{$
\def\b{\\{\hbox{}}}
\unitlength=1pt%
\def\vline(#1,#2)#3|{\leftput(#1,#2){\lline(0,1){#3}}}%
\def\hline(#1,#2)#3|{\leftput(#1,#2){\lline(1,0){#3}}}%
\ytableau{12pt}{0.4pt}{}
{\b &\b       \cr 
 \b &\b &\b   \cr
 \b &\b &\b   \cr
 \b &\b &\b &\b &\b &\b  \cr
\noalign{\vskip 3pt}
\noalign{\hbox{Fig. \sec.2. Hook length}}%
}
\vline(-82.6,27.8)37.3|
\vline(-81.0,29.4)35.8|
\vline(-79.5,30.9)34.2|
\vline(-78.0,32.4)32.6|
\vline(-76.4,34.0)31.1|
\vline(-74.8,35.6)29.6|
\vline(-73.3,37.1)28.0|
\vline(-71.8,38.6)26.4|
\vline(-70.2,40.2)24.9|
\hline(-82.8,28.0)24.9|
\hline(-81.2,29.6)23.4|
\hline(-79.7,31.1)21.8|
\hline(-78.2,32.6)20.2|
\hline(-76.6,34.2)18.7|
\hline(-75.0,35.8)17.2|
\hline(-73.5,37.3)15.6|
\hline(-72.0,38.8)14.0|
\hline(-70.4,40.4)12.5|
$
}
\setbox3=\hbox{$
\ytableau{12pt}{0.4pt}{}
{\\2 &\\1       \cr 
 \\4 &\\3 &\\1   \cr
 \\5 &\\4 &\\2   \cr
 \\9 &\\8 &\\6 &\\3 &\\2 &\\1  \cr
\noalign{\vskip 3pt}
\noalign{\hbox{Fig. \sec.3. Hook lengths}}
}$
}
$$\box1\quad\box2\quad\box3$$
}


For each box $v$ in the Ferrers diagram of a partition $\l$, or
for each box $v$ in $\l$, for short, define the 
{\it hook length} of $v$, denoted by $h_v(\l)$ or $h_v$, to be the number of 
boxes $u$ such that  $u=v$,
or $u$ lies in the same column as $v$ and above $v$, or in the 
same row as $v$ and to the right of $v$ (see Fig.~\sec.2). 
The {\it hook length multi-set} of $\l$
is the multi-set of all hook lengths of $\l$.
In Fig.~\sec.3 
the hook lengths of all boxes for the partition $\l=(6,3,3,2)$
have been written in each box. The hook length multi-set of $\l$ 
is $\{2,1,4,3,1,5,4,2,9,8,6,3,2,1\}$.

\proclaim Theorem \sec.1 [Main].
For any complex number $\beta$ we have 
$$
\prod_{m\geq 1} { (1-x^m)^{\beta-1}} 
\ =\ 
\sum_{\l\in \setP}\ \prod_{v\in\l}\bigl(1-{\beta \over h_v^2}\bigr)x.
\leqno{(\sec.6)}
$$

Identity (\sec.6) will be called ``Main Identity".
Two numerical examples are given in \S2 for verifying 
the Main Theorem. For convenience, the exponent $s$ in the Euler Product
has been replaced by $\beta-1$.
The proof of the Main Theorem is based on the Macdonald identities for
$A_l^{(a)}$. 
We have also used the properties of the bijection 
between $t$-cores and $N$-codings 
constructed by Garvan, Kim and Stanton [GKS90] (see \S5). 
From our Main Theorem we derive
new formulas about hook lengths, including the ``marked hook 
formula". We also improve a result due to Kostant (see~\S7). 
The five results we should like to single out are next stated. 
They will be further proved in Sections 2, 6, 7 and 9.

\proclaim Corollary \sec.2 [=2.4].
For any positive integers $n$ and $k$ the following expression
$$
\sum_{\l\vdash n}\ \prod_{v\in\l}\bigl(1-{k\over h_v^2}\bigr)
\leqno{(\sec.7)}
$$
is an integer.
\goodbreak

\proclaim Theorem \sec.3 [=6.10, marked hook formula].
We have
$$
\sum_{\l\vdash n} f_\l^2\sum_{v\in\l} h_v^2 = {n(3n-1)\over 2} n!,
\leqno{(\sec.8)}
$$
where $f_\l$ is the number of standard Young tableaux of shape $\l$
(see \S2.3 and \S6.4).

Theorem \sec.3 is to be compared with the following well-known formula
$$
\sum_{\l\vdash n} f_\l^2 =  n! 
\leqno{(\sec.9)}
$$

\proclaim Theorem \sec.4 [=6.9].
We have
$$
\sum_{n\geq 1} x^n \sum_{\l\vdash n} 
\bigl(\sum_{v\in\l} {1\over h_v^2}\bigr)^2
=
\prod_{m\geq 1}{1\over 1-x^m} \Bigl(
\sum_{k\geq 1}{x^k k^{-3}\over 1-x^k}+
\bigl(\sum_{k\geq 1}{x^k k^{-1}\over 1-x^k}\bigr)^2
\Bigr). \leqno{(\sec.10)}
$$

\proclaim Theorem \sec.5 [=7.2].
Let $k$ be a positive integer and $s$ be a real number 
such that $s\geq k^2-1$.  Then $(-1)^k f_k(s)>0$.

\proclaim Corollary \sec.6 [=9.2].
For any positive integer $n$ the following expression
$$
{1\over n+1}\sum_{\l\vdash n} \prod_{v\in\l} \bigl(1+{n\over h_v^2}\bigr)
\leqno{(\sec.11)}
$$
is a positive integer.

It would be interesting to find a direct proof  of 
Corollaries \sec.2, \sec.6 and Theorems \sec.3, \sec.4. 
In particular, the ``marked hook formula" suggests that 
a {\it marked} Robinson-Schensted-Knuth correspondence 
should be constructed
between
pairs of marked Young tableaux and marked permutations (see \S6.4). 
\medskip


\def\sec{2}
\section{\sec. Basic consequences and specializations} 

\noindent
\medskip \noindent
{\it \sec.1. Equivalent forms.} 
With the notations recalled above 
for partitions the right-hand side of the Main Identity can 
be written as
$$
\sum_{\l\in \setP}\ x^{|\l|}\prod_{v\in\l}\bigl(1-{\beta \over h_v^2}\bigr)
\ =\ 
\sum_{n\geq 0} x^n \sum_{\l\vdash n}\ 
\prod_{v\in\l}\bigl(1-{\beta \over h_v^2}\bigr).
\leqno{(\sec.1)}
$$
Using the identity (see [St99, p.316]) 
$$
\prod_{m\geq 1} {1\over 1-x^m} = 
\exp{\bigl(\sum_{k\geq 1} {x^k\over k(1-x^k)}\bigr)}, \leqno{(\sec.2)}
$$
the Main Identity can be written:
$$
\prod_{m\geq 1} {1\over 1-x^m}\times 
\exp{\Bigl(-\beta\sum_{k\geq 1} {x^k\over k(1-x^k)}\Bigr)} 
=
\sum_{n\geq 0} x^n \sum_{\l\vdash n}\ 
\prod_{v\in\l}\bigl(1-{\beta \over h_v^2}\bigr)
\leqno{(\sec.3)}
$$
or
$$
\exp{\Bigl((1-\beta)\sum_{k\geq 1} {x^k\over k(1-x^k)}\Bigr)}
=
\sum_{n\geq 0} x^n \sum_{\l\vdash n}\ 
\prod_{v\in\l}\bigl(1-{\beta \over h_v^2}\bigr).
\leqno{(\sec.4)}
$$
The Main Theorem has also the following equivalent generating function form,
which can be verified by comparing the coefficients of $X^mY^n$.
\proclaim Theorem \sec.1.
We have
$$
\sum_{k\geq 0}
{X^k\over 1-Y\prod_{m\geq 1}(1-x^{mk})}
=
\sum_{k\geq 0} \sum_{\l\in\setP}
{Y^k\over 1-X x^{|\l|}} \prod_{v\in\l}\Bigl( 1-{k+1\over h_v^2}\Bigr).
$$

\medskip \noindent
{\it \sec.2. Corollaries}. 
From the Main Theorem we immediately have the following results.

\proclaim Corollary \sec.2.
Let $F(\beta)$ be the function defined by
$$
F(\beta):=
\sum_{\l\in \setP}\ \prod_{v\in\l}\bigl(1-{\beta +1 \over h_v^2}\bigr)x.
$$
Then
$$
F(\beta_1 +\beta_2)=F(\beta_1) F(\beta_2). 
$$
In particular, 
$$
F(\beta) F(-\beta)=1. 
$$

\proclaim Corollary \sec.3.
For each positive integer $n$ the following expression 
$$
n!\sum_{\l\vdash n}\ \prod_{v\in\l}\bigl(1-{\beta\over h_v^2}\bigr)
$$
is a polynomial in $\beta$ with integral coefficients.

\proclaim Corollary \sec.4 [=1.2].
For any positive integers $n$ and $k$ the following expression
$$
\sum_{\l\vdash n}\ \prod_{v\in\l}\bigl(1-{k\over h_v^2}\bigr)
\leqno{(\sec.5)}
$$
is an integer.

Note that unlike  Corollary \sec.3
there is no factor $n!$ in Corollary~\sec.4. 

\medskip \noindent
{\it \sec.3. Specialization for $\beta=0$}. 
Letting $\beta=0$ in the Main Theorem yields 
the well-known generating function for partitions (see [An76, p.3]).

\proclaim Theorem \sec.5.
Let $p(n)$ be the number of partitions of $n$. Then
$$
\prod_{m\geq 1} {1\over 1-x^m}
= \sum_{n\geq 0} p(n) x^n. \leqno{(\sec.6)}
$$

Using (\sec.6) we can rewrite the Main Theorem as follows,
which is probably the good form for finding combinatorial interpretations.
It means that there exists a bijection between one weighted-partition 
(right-hand side of (\sec.7))
and sequences of $k+1$ partitions (left-hand side of (\sec.7)).
When $\beta=-k$ ($k\in\setN$), the left-hand side of the Main Identity
may be rewritten $(\sum_{n\geq 0} p(n)x^n)^{k+1}$ because of (\sec.6). 
Hence the
Main Theorem may be restated as the following corollary.

\proclaim Corollary \sec.6.
Let $k\in \setN$. Then
$$
\sum   p(n_1) p(n_2)\cdots p(n_{k+1})
\ = \ 
\sum_{\l\vdash n}\ 
\prod_{v\in\l}\bigl(1+{k \over h_v^2}\bigr),
\leqno{(\sec.7)}
$$
where the sum on the left-hand side ranges over all positive 
integer vectors $(n_1, n_2, \ldots, n_{k+1})$ such that $\sum n_i=n$.

\medskip \noindent
{\it \sec.4. Specialization for $\beta=1$}. 
The case $\beta=1$ is really trivial. Every non-empty partition $\l$
contains at least one box $v$ of hook length $h_v=1$, so that
$$
\prod_{v\in \l} \bigl(1-{1\over h_v^2} \bigr)=0.
$$
Hence
$$
\sum_{n\geq 0} x^n \sum_{\l\vdash n}\ 
\prod_{v\in\l}\bigl(1-{1\over h_v^2}\bigr)
\ =\  1. 
$$
\medskip \noindent
{\it \sec.5. Specialization for $\beta=\infty$}. 
The hook length plays an important role in Algebraic Combinatorics 
thanks to the famous hook formula
due to Frame, Robinson and Thrall [FRT54]
$$
f_\l={n!\over \prod_{v\in\l} h_v(\l)}, \leqno{(\sec.8)}
$$
where $f_\l$ is the number of standard Young tableaux of shape $\l$
(see [St99, p.376; Kn98, p.59; Kr99, Ze84, GNW79, NPS97,  RW83]).
By using the Main Theorem we re-prove the following classical result, 
which is also a consequence of the
Robinson-Schensted-Knuth 
correspondence (see, for example, [Kn98, p.49-59; St99, p.324]).

\proclaim Theorem \sec.7.
We have
$$
\sum_{\l\vdash n}\ f_\l^2 \ =\  n!  \leqno{(\sec.9)}
$$

{\it Proof}.
Put $\beta=-y/x$ in (\sec.4): 
$$
\sum_{\l\in \setP}\ \prod_{v\in\l}\Bigl(1+{y/x \over h_v^2}\Bigr)x
\ =\ 
\exp{\Bigl(\bigl(1+{y\over x}\bigr)\sum_{k\geq 1} {x^k\over k(1-x^k)}\Bigr)}. 
$$
When $x\rightarrow 0$, we get 
$$
\sum_{\l\in \setP}\ \prod_{v\in\l}{y\over h_v^2}
\ =\ 
\exp(y).  \leqno{(\sec.10)}
$$
Comparing the coefficients of $y^n$ yields 
$$
\sum_{\l\vdash n}\ \prod_{v\in\l}{1\over h_v^2}
\ =\ 
{1\over n!}  \qquad
$$
or
$$
\sum_{\l\vdash n}\ \Bigl({n!\over \prod_{v\in\l}h_v}\Bigr)^2 \ =\  n! 
\qed
$$

\medskip \noindent
{\it \sec.6. Specialization for $\beta=-1$}. 
Unlike the specializations $\beta=1, 0, \infty$ as done previously,
the case $\beta=-1$ does not relate to any classical result.
However, we get a surprising formula for 
$pp(n)$, defined to be the number of ordered pairs $\pi', \pi''$ of 
partitions such that $|\pi'|+|\pi''|=n$ { (see [BG06; CJW08])}, as stated next.

\proclaim Corollary \sec.8.
We have
$$
pp(n)
\ = \ 
\sum_{\l\vdash n}\ 
\prod_{v\in\l}\bigl(1+{1 \over h_v^2}\bigr). \leqno{(\sec.11)}
$$

Again, a direct proof of Corollary \sec.8 would be welcome.

\medskip \noindent
{\it \sec.7. Specialization for $\beta=2$}. 
There is no direct specialization for $\beta=2$. Nevertheless, by using
the Euler pentagonal theorem (formula (1.1))
we obtain another expression for the alternate sum of the pentagonal powers,
as stated in the next Proposition.

\proclaim Proposition \sec.9.
We have
$$
\sum_{\l\in \setP}\ \prod_{v\in\l}\bigl(1-{2 \over h_v^2}\bigr)x
= \sum_{k=-\infty}^{\infty}(-1)^k x^{k(3k+1)/2} .\leqno{(\sec.12)}
$$

\noindent
{\it Example \sec.1}.
We see that the coefficient of $x^4$ in (\sec.12) is $0$,
i.e.:
$$
\sum_{\l\vdash 4}\ \prod_{v\in\l}\bigl(1-{2 \over h_v^2}\bigr)
= 0.
\leqno{(\sec.13)}
$$
There are five partitions of $4$ (see. Fig.~8.1) and their hook lengths
are respectively 
$\{1,2,3,4\}, \{1, 1,2,4\}, \{1,2,2, 3\}, \{1,1,2,4\}$ and $\{1,2,3,4\}$.
We verify that (\sec.13) is true by the following calculation.
$$
\def\fh#1|{\bigl(1-{2\over #1}\bigr)}
2\fh9| \fh16| + 2\fh1| \fh16| + \fh4|\fh9| 
= 0.
$$

This raises the question: is there a direct proof of Proposition \sec.9?

\medskip \noindent
{\it \sec.8. Specialization for $\beta=25$}. 
Recall that the Ramanujan $\tau$ function is defined by (see [Se70, p.156]):
$$
\leqalignno{
x\prod_{m\geq 1} (1-x^m)^{24}&=\sum_{n\geq 1} \tau(n) x^n &{(\sec.14)}\cr
=x-24x^2+252&x^3-1472x^4+4830x^5-6048x^6+ \cdots \cr
}
$$
Putting $\beta=25$ in  the Main Theorem yields the next proposition.

\proclaim Proposition \sec.10.
We have
$$
\tau(n) \ = \ 
\sum_{\l}\ \prod_{v\in\l}\bigl(1-{25 \over h_v^2}\bigr)
\leqno{(\sec.15)}
$$
where the sum ranges over all $5$-cores of $n-1$.

\goodbreak
Notice that there is  cancellation between a box of hook length 3 and
a box of hook length 4 in each $5$-core, because
$$ 
{
\def\fh#1|{\bigl(1-{25\over #1}\bigr)}
\fh9|\fh16|=1.
}
$$

\noindent
{\it Example \sec.2}.
Take $n=6$. There are two 5-cores of 5: $(3,2)$ and $(2,2,1)$. Those two 
$5$-cores have their hook length multi-sets equal to $\{1,1,2,3,4\}$, so that
$$
{
\def\fh#1|{\bigl(1-{25\over #1}\bigr)}
\tau(6)=2\fh1|\fh1|\fh4|=-6048.
}
$$



\def\sec{3}
\section{\sec. Specialization for $\beta=4$} 
The $\beta=4$ case is very interesting. Unlike the case for $\beta=2$ in 
which Euler's pentagonal theorem is used, here the following
well-known Jacobi triple product formula is re-proved! 
(see [An76, p.21; Kn98, p.20; JS89; FH99; FK99])

\proclaim Theorem \sec.1 [Jacobi].
We have
$$
\prod_{m\geq 1} (1-x^m)^3=\sum_{m\geq 0} (-1)^m (2m+1)x^{m(m+1)/2}.
\leqno{(\sec.1)}
$$

{\it Proof}.
Put $\beta=4$ in The Main Identity:
$$
\prod_{m\geq 1} {(1-x^m)^{3}}
=
\sum_{\l\in \setP}\ \prod_{v\in\l}\bigl(1-{4 \over h_v^2}\bigr)x.
\leqno{(\sec.2)}
$$
If a partition $\l$ contains one box $v$ whose hook length is
$h_v=2$, then
$$
\prod_{v\in\l}\bigl(1-{4 \over h_v^2}\bigr)x =0.\leqno{(\sec.3)}
$$
Otherwise $\l$ must be a {\it staircase partition} 
$$\Delta_m:=(m, m-1, \ldots, 3,2,1).$$ 
We have (see Fig. \sec.1 and \sec.2 for an example):
$$
\leqalignno{
\prod_{v\in \Delta_m} \bigl(1-{4\over h_v^2}\bigr) \
&= 
\Bigl({ (2m-1)^2-4\over (2m-1)^2} \Bigr)^1
\cdots
\Bigl({ 5\over 9 }\Bigr)^{m-1}
\Bigl({ -3\over 1} \Bigr)^{m}\cr
& = (-1)^m (2m+1). &{(\sec.4)}\cr
}
$$

\medskip
\setbox1=\hbox{$
\def\b{\\{\hbox{}}}
\ytableau{18pt}{0.4pt}{}
{\\1        \cr 
 \\3 &\\1    \cr
 \\5 &\\3 & \\1   \cr
 \\7 &\\5 &\\3  &\\1  \cr
\noalign{\vskip 3pt}
\noalign{\hbox{\kern -7mm Fig. \sec.1. Hook lengths $h_v$ in $\Delta_4$}}
}$
}
\setbox2=\hbox{$
\ytableau{18pt}{0.4pt}{}
{\\{-3}       \cr 
 \\{5\over9}  &\\{-3}   \cr
 \\{21\over 25} &\\{5\over9} &\\{-3}   \cr
 \\{45\over 49} &\\{21\over 25} &\\{5\over 9} &\\{-3}  \cr
\noalign{\vskip 3pt}
\noalign{\hbox{\kern -2mm Fig. \sec.2. $(1-4/h_v^2)$ in $\Delta_4$}}
}$
}
\midinsert
$$
\box1\qquad\qquad\box2
$$
\endinsert

\noindent
As $|\Delta_m|=m(m+1)/2$, we conclude that
$$
\prod_{m\geq 1} {(1-x^m)^{3}}
=\sum_{m\geq 0} \prod_{v\in \Delta_m}\bigl(1-{4 \over h_v^2}\bigr)x
= \sum_{m\geq 0} (-1)^m (2m+1)x^{m(m+1)/2}.\qed
$$


\def\sec{4}
\section{\sec. Specialization for $\beta=9$} 
Recall that a partition 
$\l$ is a {\it $t$-core}, if $\l$ has no hook length equal to $t$. 
[GKS90; St99, p.468]. Hence, if $\l$ is not a 3-core, 
$$
\prod_{v\in\l}\bigl(1-{9 \over h_v^2}\bigr)x =0. \leqno{(\sec.1)}
$$
By the Main Theorem
$$
\prod_{m\geq 1} { (1-x^m)^8} \ = \ 
\sum_{\l}\ \prod_{v\in\l}\bigl(1-{9 \over h_v^2}\bigr)x
\leqno{(\sec.2)}
$$
where the sum ranges over all $3$-cores.
\smallskip

\proclaim Theorem \sec.1.
We have
$$
\leqalignno{
\prod_{k\geq 1} (1-q^k)^8
&=
\sum_{k, m\geq 0} \Bigl(
{1\over 2}   (3k+1)(3m+1)(3k+3m+2)
q^{k^2+k+m^2+m+km} \cr
&\ 
-{1\over 2}   (3k+2)(3m+2)(3k+3m+4)
q^{k^2+k+m^2+m+(k+1)(m+1)} \Bigr). \cr
}
$$

{\it Proof}.
We need characterize all $3$-cores.
In fact, a partition is a $3$-core if and only if
it has one of the forms described in Fig.~\sec.1 (type A) and
Fig.~\sec.2 (type B). 
Let $\Delta_k=(k,\ldots, 3,2,1)$ be a staircase partition.
Define 
$$\Delta^2_k=(k,k, \ldots, 3,3,2,2,1,1)$$ 
and ${\Delta^{2}_k}'$ be the transposition of $\Delta^2_k$.
Then type A (resp. type B) is made of a partition of form $\Delta^2_k$, 
a partition of form ${\Delta^{2}_m}'$ and a rectangle of form $k\times m$
(resp. of form $(k+1)\times(m+1)$).
We have (see Fig. \sec.3 and Fig. \sec.4)
%
{
\setbox1=\hbox{$
\def\b{\\{\hbox{}}}
\ytableau{12pt}{0.4pt}{}
{
 \b      \cr 
 \b        \cr 
 \b &\b       \cr 
 \b &\b       \cr 
 \b &\b &\b   \cr
 \b &\b &\b   \cr
 \b &\b &\b &\b  \cr
 \b &\b &\b &\b  \cr
 \\{\hbox{14}} &\\{\hbox{11}}&\\{\hbox{8}}&\\{\hbox{5}}&\b &\b  \cr
 \\{\hbox{17}} &\\{\hbox{14}}&\\{\hbox{11}}&\\{\hbox{8}}&\b &\b & \b &\b\cr
 \\{\hbox{20}} &\\{\hbox{17}}&\\{\hbox{14}}&\\{\hbox{11}}&\b &\b & 
         \b &\b &\b &\b\cr
\noalign{\vskip 3pt}
\noalign{\hbox{\quad Fig. \sec.1. Type A $3$-core}}
}$
}
\setbox2=\hbox{$
\def\b{\\{\hbox{}}}
\ytableau{12pt}{0.4pt}{}
{
 \b      \cr 
 \b        \cr 
 \b &\b       \cr 
 \b &\b       \cr 
 \b &\b &\b   \cr
 \b &\b &\b   \cr
 \b &\b &\b &\b  \cr
 \b &\b &\b &\b  \cr
 \\{\hbox{13}}&\\{\hbox{10}}&\\{\hbox{7}}&\\{\hbox{4}}&\\{\hbox{1}}&\cr
 \\{\hbox{16}}&\\{\hbox{13}}&\\{\hbox{10}}&\\{\hbox{7}}&\\{\hbox{4}}&\b &\b\cr
 \\{\hbox{19}} &\\{\hbox{16}}&\\{\hbox{13}}&\\{\hbox{10}}&\\{\hbox{7}}
     & \b &\b &\b &\b\cr
 \\{\hbox{22}} &\\{\hbox{19}} &\\{\hbox{16}}&\\{\hbox{13}}&\\{\hbox{10}}&\b 
     &\b & \b &\b &\b &\b\cr
\noalign{\vskip 3pt}
\noalign{\hbox{\quad Fig. \sec.2. Type B $3$-core}}
}$
}

\midinsert
\vskip -10pt
$$
\box1\qquad\box2
$$
\vskip -10pt
\endinsert
}

\setbox1=\hbox{$
\def\b{\\{\hbox{}}}
\ytableau{18pt}{0.4pt}{}
{\\2 & \\1       \cr  
 \\5 &\\4 &\\2 &\\1   \cr
 \\8 &\\7 &\\5 &\\4 & \\2 &\\1  \cr
\noalign{\vskip 3pt}
\noalign{\hbox{\kern -2mm Fig. \sec.3. Hook lengths $h_v$ in ${\Delta^2_3}'$}}
}$
}
\setbox2=\hbox{$
\ytableau{18pt}{0.4pt}{}
{\\{-{5\over 4}}  & \\{-8}     \cr 
 \\{16\over 25} & \\{7\over 16} &\\{-{5\over4}}  &\\{-8}   \cr
 \\{55\over 64}& \\{40\over 49} &\\{16\over 25} & \\{7\over 16} &\\{-{5\over4}}  &\\{-8}   \cr
\noalign{\vskip 3pt}
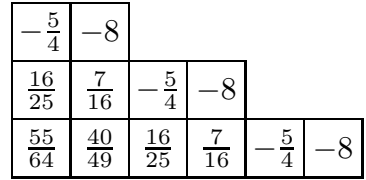
\noalign{\hbox{\kern 0mm Fig. \sec.4. $(1-9/h_v^2)$ in ${\Delta^2_3}'$}}
}$
}
\midinsert
\vskip -10pt
$$
\box1\qquad\qquad\box2
$$
\vskip -10pt
\endinsert

$$
\leqalignno{
\prod_{v\in\Delta^2_k} \Bigl(1- {9\over h_v^2}\Bigr)
&=
(-8)^k \Bigl(-{5\over 4}\Bigr)^k \Bigl({7\over 16}\Bigr)^{k-1} 
\Bigl({16\over 25}\Bigr)^{k-1} \cdots \cr
&\qquad \times \Bigl({(3k-2)^2-9 \over (3k-2)^2}\Bigr)^1
\Bigl({(3k-1)^2-9\over (3k-1)^2}\Bigr)^1 \cr
&= (3k+1)(3k+2)/2. &(\sec.3)\cr
}
$$
The product of $1-9/h^2_v$ for all boxes $v$ in the rectangle of 
type A (Fig.~\sec.1) 
is:
$$
\leqalignno{
\prod_{v\in A(k,m)} \Bigl(1- {9\over h_v^2}\Bigr)
&=
\prod_{j=1}^m \prod_{i=1}^k {(3i+3j-1)^2-9 \over (3i+3j-1)^2} \cr
&=\prod_{j=1}^m {(3j-1) (3j+3k+2)\over (3j+2)(3j+3k-1) }\cr
&={2 (3k+3m+2)\over (3k+2)(3m+2) }. &(\sec.4)\cr
}
$$
The product of $1-9/h^2_v$ for all boxes $v$ in the rectangle of 
type B (Fig.~\sec.2)
is:
$$
\leqalignno{
\prod_{v\in B(k,m)} \Bigl(1- {9\over h_v^2}\Bigr)
&=
\prod_{j=1}^{m+1} \prod_{i=1}^{k+1} {(3i+3j-5)^2-9 \over (3i+3j-5)^2} \cr
&={-2 (3k+3m+4)\over (3k+1)(3m+1) }. &(\sec.5)\cr
}
$$
Combining equations (\sec.3), (\sec.4) and (\sec.5) yields Theorem \sec.1.
\qed


\def\sec{5}
\section{\sec. Proof of the Main Theorem} 
In this section we always suppose that $t=2t'+1$ is an odd positive integer. 
\medskip \noindent
{\it \sec.1. Fundamental properties of $t$-cores and $V$-codings}. 
Recall that a partition $\l$ is a $t$-core if the {\it hook length multi-set} 
of $\l$ does not contain the integer $t$.
It is known that the hook length multi-set of each $t$-core 
does not contain any {\it multiple} of $t$ 
[Kn98. p.69, p.612; St99, p.468].

\medskip

{\it Definition \sec.1}.
Each vector of integers $(v_0, v_1, \ldots, v_{t-1})\in \setZ^t$ is called 
{\it $V$-coding} if the following conditions hold:

(i) $v_i\equiv i (\mod t)$ for $0\leq i\leq t-1$; 

(ii) $v_0+v_1+\cdots + v_{t-1}=0$.

\smallskip
The $V$-coding can be identified with the {\it set} 
$\{ v_0, v_1,  \ldots, v_{t-1}\}$ thanks to condition (i).
In this section we present a bijection between $t$-cores and $V$-codings,
that constitutes the crucial step in the proof of the Main Theorem.
\medskip

\proclaim Theorem \sec.1.
There is a bijection $\phi_V: \l\mapsto (v_0, v_1, \ldots, v_{t-1})$ 
which maps each $t$-core onto a $V$-coding such that
$$
|\l|={1\over 2t}(v_0^2+v_1^2+\cdots +v_{t-1}^2) - {t^2-1\over 24}
\leqno{(\sec.1)}
$$
and
$$
\prod_{v\in\l} \Bigl(1-{t^2\over h_v^2}\Bigr)
={(-1)^{t'}\over 1!\cdot 2!\cdot 3!\cdots (t-1)!} 
\prod_{0\leq i<j\leq t-1} (v_i-v_j).
\leqno{(\sec.2)}
$$

We will describe the bijection $\phi_V$ and
prove the two equalities (\sec.1) and (\sec.2) in \S\sec.2, \S\sec.3 and
\S\sec.4 respectively.
An example is given after
the construction of the bijection $\phi_V$.
\goodbreak\medskip\noindent
{\it \sec.2. The bijection $\phi_V$ and an example}. 
Each finite set of integers $A=\{a_1, a_2, \ldots, a_n\}$ is said to be 
{\it $t$-compact} if the following conditions hold:

(i) $-1, -2, \ldots, -t\in A$;

(ii) for each $a\in A$ such that 
$a\not=-1,-2, \ldots, -t$, we have $a\geq 1$ and $a \not\equiv 0 \mod t$; 

(iii) let $b>a\geq 1$ be two integers such that $a \equiv b \mod t$. 
If $b\in A$,  then $a\in A$.
\smallskip

Let $A$ be a $t$-compact set. An element $a\in A$ is said to be 
{\it $t$-maximal}
if $b\not\in A$ for every $b>a$ such that $a\equiv b\mod t$.
The set of $t$-maximal letters of $A$ is denoted by
$\tmax(A)$. 
Let $\l$ be a $t$-core.
The {\it $H$-set} of the $t$-core $\l$ is defined to be 
$$
H(\l)=\{h_v \mid \hbox{$v$ is a box in the leftmost 
column of $\l$}\}\cup\{-1,-2, \ldots -t\}.
$$

\medskip

\proclaim Lemma \sec.2.
For each $t$-core $\l$ its $H$-set 
$H(\l)$ is a $t$-compact set.

{\it Proof}. Let $c=tk+r$ ($k\geq 1, 0\leq r\leq t-1$) be an element in $H(\l)$ 
and $a$ be the maximal letter in $H(\l)$ such that $a<t(k-1)+r$. 
We must show that $t(k-1)+r$ is also in $H(\l)$. 
If it were not the case, let 
$z>t(k-1)+r, y_1, y_2, \ldots, y_d$ be the hook lengths as shown in 
Fig. \sec.1, 
where only the relevant horizontal section of the 
partition diagram has been represented.
We have
$y_1=c-a-1 \geq tk+r - t(k-1)-r  = t$
and 
$y_d=c-z+1 \leq tk+r - t(k-1)-r = t$;
so that there is one hook $y_i=t$. 
This is a contradiction since $\l$
is supposed to be a $t$-core. 
\qed

\medskip 
{
\setbox1=\hbox{$
\def\b{\\{\hbox{}}}
\ytableau{14pt}{0.4pt}{}
{
 \\a  &\b  &\b  &\b  &\b \cr
 \\z   &\b  &\b  &\b  &\b  &\b  &\b &\b  &\b  \cr
 \b   &\b  &\b  &\b  &\b  &\b  &\b &\b  &\b  &\b \cr
 \b   &\b  &\b  &\b  &\b  &\b  &\b &\b  &\b  &\b \cr
 \\c  &\b  &\b  &\b  &\b &\\{y_1}  &\\{y_2} &\\{\cdots}  &\\{y_d} 
     &\b &\b  &\b  \cr
\noalign{\vskip 3pt}
\noalign{\hbox{Fig. \sec.1. Hook length and $t$-compact set}}
}$
}
\midinsert
$$
\box1
$$
\endinsert
}

{\it Construction of $\phi_V$}. Let $\l$ be a $t$-core and $H(\l)$ be
its $H$-set.
The {\it $U$-coding} of $\l$ is defined to be the set
$U:=\tmax(H(\l))$, which can be identified with the vector
$(u_0, u_1, \ldots, u_{t-1})$ such that 
$u_0=-t$, $u_i>-t$ and $u_i\equiv i \mod t$ for $1\leq i\leq t-1$.
In general,
$$S:=u_0 + u_1 +\cdots + u_{t-1} \not=0.$$ 
The integer $S$ is a multiple of $t$ because
$$
S=\sum u_i = \sum (tk_i +i) = t\sum k_i + t(t-1)/2 
$$
(remember that $t=2t'+1$ is an odd integer). 
The $V$-coding $\phi_V(\l)$ is the set $V$ obtained from $U$ by 
the following {\it normalization}: 
$$
\phi_V(\l)=V:=\{u-S/t \ : \ u\in U\}. \leqno{(\sec.3)}
$$
In fact, we can prove that $S/t=\ell(\l)-t'-1$  (see (\sec.8)).
The set $V$ can be identified with a vector $V$-coding because 
$$\sum v_i=\sum (u_i -S/t)=\sum u_i - S=0.$$
\medskip

\noindent
{\it Example \sec.1}.
Consider the $5$-core 
$$\l=(14,10,6,6,4,4,4,2,2,2).$$ 
The $H$-set of $\l$ (see Fig. \sec.2)
$$
H(\l)=\{23,18,13,12,9,8,7,4,3,2, -1, -2, -3, -4, -5\}
$$
is $5$-compact. 
The $U$-coding of $\l$ is $U=\max_5(H(\l))=\{23,12,9,-4,-5\}$,
or in vector form  
$$(u_0,u_1,u_2,u_3,u_4)=(-5, -4, 12, 23, 9).$$
As $S=\sum u_i=35$,  the $V$-coding is given by
$$V=\{ -5-7, -4-7, 12-7, 23-7, 9-7 \}
=\{ -12, -11, 5, 16, 2 \},$$
or in vector form
$$
\phi_V(\l)=(v_0, v_1, v_2, v_3, v_4)=(5, 16, 2, -12, -11).
$$



\long\def\maplebegin#1\mapleend{}

\maplebegin

# --------------- begin maple ----------------------

# Copy the following text  to "makefig.mpl"
# then in maple > read("makefig.mpl");
# it will create a file "z_fig_by_maple.tex"

#\unitlength=1pt

Hu:= 16.4; # height quantities
Lu:= Hu; # large unity

X0:=0; Y0:=40; # origin position

File:=fopen("z_fig_by_maple.tex", WRITE);

hdash:=proc(x,y,len)
local i, d,sp, xx;
	d:=2;
	sp:=2;
	for xx from x*Lu+X0 to x*Lu+X0+Lu*len-d by d+sp do
	fprintf(File, "\\hline(
	od:
end;

vdash:=proc(x,y,len)
local i, d,sp, yy;
	d:=2;
	sp:=2;
	for yy from y*Hu+Y0 to y*Hu+Y0+Hu*len-d by d+sp do
		fprintf(File, "\\vline(
	od:
end;

text:=proc(x,y,t)
	fprintf(File, "\\centerput(
end;

cercle:=proc(x,y) text(x,y, "\\cerclechar"); end;
# the partition

# h-line

cercle(0,14);    cercle(-1,14);   
cercle(0,13);    cercle(-1,13);
cercle(0,4);     cercle(3,4);
cercle(0,3);     cercle(5,3);
cercle(0,0);     cercle(13,0);

text(-1,14, "{\\it 0}");text(0,14, "-5");
text(-1,13, "{\\it 1}");text(0,13, "-4");
text(-1,12, "{\\it 2}");text(0,12, "-3");
text(-1,11, "{\\it 3}");text(0,11, "-2");
text(-1,10, "{\\it 4}");text(0,10, "-1");
text(-1,9, "{\\it 0}");

hdash(-1,0,1);
hdash(-1,15,2);
vdash(-1,0,15);
vdash(0,10,5);
vdash(1,10,5);

# region separ
hhdash:=proc(x,y) hdash(x,y+0.03,1); end;
HHdash:=proc(XL, dx) local x; for x in XL do hhdash(x, x+dx); od; end;
vvdash:=proc(x,y) vdash(x,y+0.03,1); end;
VVdash:=proc(XL, dx) local x; for x in XL do vvdash(x, x+dx); od; end;

HHdash([-1,0,1,2,3], 11);
VVdash([   0,1,2,3], 10);

HHdash([-1,0,1,2,3,4,5,6], 6);
VVdash([   0,1,2,3,4,5,6], 5);

HHdash([-1,0,1,2,3,4,5,6,7,8], 1);
VVdash([   0,1,2,3,4,5,6,7,8], 0);

HHdash([4,5,6,7,8,9,10,11], -4);
VVdash([  5,6,7,8,9,10,11], -5);

HHdash([9,10,11,12,13], -9);
VVdash([  10,11,12,13], -10);

text(2,14.2, "$r=-2$");
text(4,12.2, "$r=-1$");
text(6,9.8, "$r=0$");
text(8,6.8, "$r=1$");
text(11,3.8, "$r=2$");
text(12.5,1.5, "$r=3$");

fclose(File);
# -------------------- end maple -------------------------
\mapleend


{
\unitlength=1pt
\font\cerclefont=cmsy10 at 14pt
\def\cerclechar{\hbox{\cerclefont\char'015}}
\def\vline(#1,#2)#3|{\leftput(#1,#2){\lline(0,1){#3}}}
\def\hline(#1,#2)#3|{\leftput(#1,#2){\lline(1,0){#3}}}
\newbox\boxhook
\setbox\boxhook=\vbox{
\vskip 26mm
$
\def\b{\\{\hbox{}}}
\def\r#1{\\{{\it #1}}}
\hskip 48pt \ytableau{16pt}{0.4pt}{}
{
 \\2 &\r2       \cr 
 \\3 &\r3       \cr 
 \\4 &\r4       \cr 
 \\7 &\r0 &\r1 &\r2  \cr
 \\8 &\b &\b &\r3  \cr
 \\9 &\b &\b &\r4  \cr
 \\{12} &\b &\b &\r0 &\r1 &\r2 \cr
 \\{13} &\b &\b &\b &\b &\r3 \cr
 \\{18} &\b &\b &\b &\b &\r4 &\r0 &\r1 &\r2 &\r3 \cr
 \\{23} &\b &\b &\b &\b &\b &\b &\b
        &\b &\r4 &\r0 &\r1 &\r2 &\r3 \cr
\noalign{\vskip 3pt}
\noalign{\hbox{\quad Fig. \sec.2. $U$-coding and $N$-coding of $t$-core}}
\centerput(0,270){\cerclechar}
\centerput(-16,270){\cerclechar}
\centerput(0,254){\cerclechar}
\centerput(-16,254){\cerclechar}
\centerput(0,106){\cerclechar}
\centerput(50,106){\cerclechar}
\centerput(0,90){\cerclechar}
\centerput(82,90){\cerclechar}
\centerput(0,40){\cerclechar}
\centerput(214,40){\cerclechar}
\centerput(-16,270){{\it 0}}
\centerput(0,270){-5}
\centerput(-16,254){{\it 1}}
\centerput(0,254){-4}
\centerput(-16,237){{\it 2}}
\centerput(0,237){-3}
\centerput(-16,221){{\it 3}}
\centerput(0,221){-2}
\centerput(-16,204){{\it 4}}
\centerput(0,204){-1}
\centerput(-16,188){{\it 0}}
\hline(-24,36)2|
\hline(-20,36)2|
\hline(-16,36)2|
\hline(-12,36)2|
\hline(-24,282)2|
\hline(-20,282)2|
\hline(-16,282)2|
\hline(-12,282)2|
\hline(-8,282)2|
\hline(-4,282)2|
\hline(0,282)2|
\hline(4,282)2|
\vline(-24,36)2|
\vline(-24,40)2|
\vline(-24,44)2|
\vline(-24,48)2|
\vline(-24,52)2|
\vline(-24,56)2|
\vline(-24,60)2|
\vline(-24,64)2|
\vline(-24,68)2|
\vline(-24,72)2|
\vline(-24,76)2|
\vline(-24,80)2|
\vline(-24,84)2|
\vline(-24,88)2|
\vline(-24,92)2|
\vline(-24,96)2|
\vline(-24,100)2|
\vline(-24,104)2|
\vline(-24,108)2|
\vline(-24,112)2|
\vline(-24,116)2|
\vline(-24,120)2|
\vline(-24,124)2|
\vline(-24,128)2|
\vline(-24,132)2|
\vline(-24,136)2|
\vline(-24,140)2|
\vline(-24,144)2|
\vline(-24,148)2|
\vline(-24,152)2|
\vline(-24,156)2|
\vline(-24,160)2|
\vline(-24,164)2|
\vline(-24,168)2|
\vline(-24,172)2|
\vline(-24,176)2|
\vline(-24,180)2|
\vline(-24,184)2|
\vline(-24,188)2|
\vline(-24,192)2|
\vline(-24,196)2|
\vline(-24,200)2|
\vline(-24,204)2|
\vline(-24,208)2|
\vline(-24,212)2|
\vline(-24,216)2|
\vline(-24,220)2|
\vline(-24,224)2|
\vline(-24,228)2|
\vline(-24,232)2|
\vline(-24,236)2|
\vline(-24,240)2|
\vline(-24,244)2|
\vline(-24,248)2|
\vline(-24,252)2|
\vline(-24,256)2|
\vline(-24,260)2|
\vline(-24,264)2|
\vline(-24,268)2|
\vline(-24,272)2|
\vline(-24,276)2|
\vline(-24,280)2|
\vline(-8,200)2|
\vline(-8,204)2|
\vline(-8,208)2|
\vline(-8,212)2|
\vline(-8,216)2|
\vline(-8,220)2|
\vline(-8,224)2|
\vline(-8,228)2|
\vline(-8,232)2|
\vline(-8,236)2|
\vline(-8,240)2|
\vline(-8,244)2|
\vline(-8,248)2|
\vline(-8,252)2|
\vline(-8,256)2|
\vline(-8,260)2|
\vline(-8,264)2|
\vline(-8,268)2|
\vline(-8,272)2|
\vline(-8,276)2|
\vline(-8,280)2|
\vline(9,200)2|
\vline(9,204)2|
\vline(9,208)2|
\vline(9,212)2|
\vline(9,216)2|
\vline(9,220)2|
\vline(9,224)2|
\vline(9,228)2|
\vline(9,232)2|
\vline(9,236)2|
\vline(9,240)2|
\vline(9,244)2|
\vline(9,248)2|
\vline(9,252)2|
\vline(9,256)2|
\vline(9,260)2|
\vline(9,264)2|
\vline(9,268)2|
\vline(9,272)2|
\vline(9,276)2|
\vline(9,280)2|
\hline(-24,201)2|
\hline(-20,201)2|
\hline(-16,201)2|
\hline(-12,201)2|
\hline(-8,217)2|
\hline(-4,217)2|
\hline(0,217)2|
\hline(4,217)2|
\hline(9,234)2|
\hline(13,234)2|
\hline(17,234)2|
\hline(21,234)2|
\hline(25,250)2|
\hline(29,250)2|
\hline(33,250)2|
\hline(37,250)2|
\hline(42,267)2|
\hline(46,267)2|
\hline(50,267)2|
\hline(54,267)2|
\vline(-8,201)2|
\vline(-8,205)2|
\vline(-8,209)2|
\vline(-8,213)2|
\vline(9,217)2|
\vline(9,221)2|
\vline(9,225)2|
\vline(9,229)2|
\vline(25,234)2|
\vline(25,238)2|
\vline(25,242)2|
\vline(25,246)2|
\vline(42,250)2|
\vline(42,254)2|
\vline(42,258)2|
\vline(42,262)2|
\hline(-24,119)2|
\hline(-20,119)2|
\hline(-16,119)2|
\hline(-12,119)2|
\hline(-8,135)2|
\hline(-4,135)2|
\hline(0,135)2|
\hline(4,135)2|
\hline(9,152)2|
\hline(13,152)2|
\hline(17,152)2|
\hline(21,152)2|
\hline(25,168)2|
\hline(29,168)2|
\hline(33,168)2|
\hline(37,168)2|
\hline(42,185)2|
\hline(46,185)2|
\hline(50,185)2|
\hline(54,185)2|
\hline(58,201)2|
\hline(62,201)2|
\hline(66,201)2|
\hline(70,201)2|
\hline(74,217)2|
\hline(78,217)2|
\hline(82,217)2|
\hline(86,217)2|
\hline(91,234)2|
\hline(95,234)2|
\hline(99,234)2|
\hline(103,234)2|
\vline(-8,119)2|
\vline(-8,123)2|
\vline(-8,127)2|
\vline(-8,131)2|
\vline(9,135)2|
\vline(9,139)2|
\vline(9,143)2|
\vline(9,147)2|
\vline(25,152)2|
\vline(25,156)2|
\vline(25,160)2|
\vline(25,164)2|
\vline(42,168)2|
\vline(42,172)2|
\vline(42,176)2|
\vline(42,180)2|
\vline(58,185)2|
\vline(58,189)2|
\vline(58,193)2|
\vline(58,197)2|
\vline(74,201)2|
\vline(74,205)2|
\vline(74,209)2|
\vline(74,213)2|
\vline(91,217)2|
\vline(91,221)2|
\vline(91,225)2|
\vline(91,229)2|
\hline(-24,37)2|
\hline(-20,37)2|
\hline(-16,37)2|
\hline(-12,37)2|
\hline(-8,53)2|
\hline(-4,53)2|
\hline(0,53)2|
\hline(4,53)2|
\hline(9,70)2|
\hline(13,70)2|
\hline(17,70)2|
\hline(21,70)2|
\hline(25,86)2|
\hline(29,86)2|
\hline(33,86)2|
\hline(37,86)2|
\hline(42,103)2|
\hline(46,103)2|
\hline(50,103)2|
\hline(54,103)2|
\hline(58,119)2|
\hline(62,119)2|
\hline(66,119)2|
\hline(70,119)2|
\hline(74,135)2|
\hline(78,135)2|
\hline(82,135)2|
\hline(86,135)2|
\hline(91,152)2|
\hline(95,152)2|
\hline(99,152)2|
\hline(103,152)2|
\hline(107,168)2|
\hline(111,168)2|
\hline(115,168)2|
\hline(119,168)2|
\hline(124,185)2|
\hline(128,185)2|
\hline(132,185)2|
\hline(136,185)2|
\vline(-8,37)2|
\vline(-8,41)2|
\vline(-8,45)2|
\vline(-8,49)2|
\vline(9,53)2|
\vline(9,57)2|
\vline(9,61)2|
\vline(9,65)2|
\vline(25,70)2|
\vline(25,74)2|
\vline(25,78)2|
\vline(25,82)2|
\vline(42,86)2|
\vline(42,90)2|
\vline(42,94)2|
\vline(42,98)2|
\vline(58,103)2|
\vline(58,107)2|
\vline(58,111)2|
\vline(58,115)2|
\vline(74,119)2|
\vline(74,123)2|
\vline(74,127)2|
\vline(74,131)2|
\vline(91,135)2|
\vline(91,139)2|
\vline(91,143)2|
\vline(91,147)2|
\vline(107,152)2|
\vline(107,156)2|
\vline(107,160)2|
\vline(107,164)2|
\vline(124,168)2|
\vline(124,172)2|
\vline(124,176)2|
\vline(124,180)2|
\hline(58,37)2|
\hline(62,37)2|
\hline(66,37)2|
\hline(70,37)2|
\hline(74,53)2|
\hline(78,53)2|
\hline(82,53)2|
\hline(86,53)2|
\hline(91,70)2|
\hline(95,70)2|
\hline(99,70)2|
\hline(103,70)2|
\hline(107,86)2|
\hline(111,86)2|
\hline(115,86)2|
\hline(119,86)2|
\hline(124,103)2|
\hline(128,103)2|
\hline(132,103)2|
\hline(136,103)2|
\hline(140,119)2|
\hline(144,119)2|
\hline(148,119)2|
\hline(152,119)2|
\hline(156,135)2|
\hline(160,135)2|
\hline(164,135)2|
\hline(168,135)2|
\hline(173,152)2|
\hline(177,152)2|
\hline(181,152)2|
\hline(185,152)2|
\vline(74,37)2|
\vline(74,41)2|
\vline(74,45)2|
\vline(74,49)2|
\vline(91,53)2|
\vline(91,57)2|
\vline(91,61)2|
\vline(91,65)2|
\vline(107,70)2|
\vline(107,74)2|
\vline(107,78)2|
\vline(107,82)2|
\vline(124,86)2|
\vline(124,90)2|
\vline(124,94)2|
\vline(124,98)2|
\vline(140,103)2|
\vline(140,107)2|
\vline(140,111)2|
\vline(140,115)2|
\vline(156,119)2|
\vline(156,123)2|
\vline(156,127)2|
\vline(156,131)2|
\vline(173,135)2|
\vline(173,139)2|
\vline(173,143)2|
\vline(173,147)2|
\hline(140,37)2|
\hline(144,37)2|
\hline(148,37)2|
\hline(152,37)2|
\hline(156,53)2|
\hline(160,53)2|
\hline(164,53)2|
\hline(168,53)2|
\hline(173,70)2|
\hline(177,70)2|
\hline(181,70)2|
\hline(185,70)2|
\hline(189,86)2|
\hline(193,86)2|
\hline(197,86)2|
\hline(201,86)2|
\hline(206,103)2|
\hline(210,103)2|
\hline(214,103)2|
\hline(218,103)2|
\vline(156,37)2|
\vline(156,41)2|
\vline(156,45)2|
\vline(156,49)2|
\vline(173,53)2|
\vline(173,57)2|
\vline(173,61)2|
\vline(173,65)2|
\vline(189,70)2|
\vline(189,74)2|
\vline(189,78)2|
\vline(189,82)2|
\vline(206,86)2|
\vline(206,90)2|
\vline(206,94)2|
\vline(206,98)2|
\centerput(33,273){$r=-2$}
\centerput(66,241){$r=-1$}
\centerput(99,201){$r=0$}
\centerput(132,152){$r=1$}
\centerput(181,103){$r=2$}
\centerput(205,65){$r=3$}
}$
}
\midinsert
$$
\box\boxhook
$$
\vskip -20pt
\endinsert
}
\noindent
We have
$$
\leqalignno{
|\l|&={1\over 2t}(v_0^2+v_1^2+\cdots +v_{t-1}^2) - {t^2-1\over 24}\cr
&= {1\over 2\cdot 5} ( 5^2+ 16^2+ 2^2+(-12)^2+(-11)^2) - {5^2-1\over 24}=54.\cr
}
$$
and
$$
\leqalignno{
\prod_{v\in\l} \Bigl(1-{5^2\over h_v^2}\Bigr)
&={1\over 1!\cdot 2!\cdot 3!\cdots (t-1)!} 
\prod_{0\leq i<j\leq t-1} (v_i-v_j)\cr
&=(-11)(3)(17)(16)\cdot(14)(28)(27)\cdot(14)(13)\cdot(-1)/288\cr
&=60035976.\cr
}
$$
Notice that, as expected, the above two numbers are positive integers.
\goodbreak\noindent
{\it \sec.3. Proof of (\sec.1)}. 
A vector of integers $(n_0, n_1, \ldots, n_{t-1})\in \setZ^t$ is said to be an 
{\it $N$-coding} if $n_0+n_1+\cdots+n_{t-1}=0$. 
Garvan, Kim and Stanton have defined a bijection $\phi_N$ between
$N$-codings and $t$-cores. We now recall its definition using their own words 
[GKS90,p.3] (see also [BG06]).

Let $\l$ be a $t$-core. Define the vector $(n_0, \ldots, n_{t-1})=\phi_N(\l)$
in the following way. Label a box in the $i$-th row and $j$-column of $\l$
by $j-i \mod t$. 
We also label the boxes in column 0 (in dotted lines in Fig. \sec.2) in
the same way, and call the resulting diagram the {\it extended $t$-residue
diagram}. A box is called {\it exposed} if it is at the end of a row of the
extended $t$-residue diagram. 
The set of boxes $(i,j)$ satisfying 
$t(r-1)\leq j-i < tr$ of the extended $t$-residue diagram of $\l$ is 
called {\it region} and numbered $r$. 
In Fig. \sec.2 the regions have been bordered by 
dotted lines.  
We now define $n_i$ to be the 
maximum region $r$ which contains an exposed box labeled $i$.
\medskip

In Fig. \sec.2 the labels of all boxes lying on the maximal border strip
(but the leftmost one) have been written in italic. This includes all the
exposed boxes: 3,3,3,2,4,3,2,4,3,2,4,3,2,1,0, when reading from bottom to top.
We have 
$(n_0, n_1, n_2, n_3, n_4)=(-2, -2, 1, 3, 0)$.

\proclaim Theorem \sec.3 [Garvan-Kim-Stanton].
The bijection 
$$\phi_N : \l \mapsto (n_0, n_1, \ldots, n_{t-1})$$
has the following property:
$$
|\l| = {t\over 2} \sum_{i=0}^{t-1} n_i^2 + \sum_{i=0}^{t-1} in_i. 
\leqno{(\sec.4)}
$$

Let $t'=(t-1)/2$ and let
$$\phi_V^N : (n_0, n_1, \ldots, n_{t-1}) \mapsto
(v_0, v_1, \ldots, v_{t-1})$$
be the bijection that maps each $N$-coding onto the $V$-coding defined by
$$
v_i=\cases
{
t n_{i+t'} +i & if $0\leq i\leq t'$; \cr
t n_{i-t'-1} +i-t & if $t'+1\leq i\leq t-1$ \cr
}
\leqno{(\sec.5)}
$$
or in set form
$$
\{v_i \mid 0\leq i \leq t-1\}= \{ t n_i +i -t' \mid 0\leq i \leq t-1\}.
\leqno{(\sec.6)}
$$
The bijective property is easy to verify.
More essentially, the bijection $\phi_V$ defined in \S5.2 is the composition 
product of the two previous bijections as is now shown.

\proclaim Lemma \sec.4.
We have $\phi_V = \phi_V^N\circ\phi_N$.

{\it Proof}. 
Let $(v_0, \ldots, v_{t-1})=\phi_V(\l)$,
$(n_0, \ldots, n_{t-1})=\phi_N(\l)$ and
$$(v_0', \ldots, v_{t-1}')=\phi_V^N(n_0, \ldots, n_{t-1}).$$
We need prove that $v_i=v_i'$.
The number $n_i$ in the $N$-coding is
defined to be the maximum region $r$ which contains an exposed box 
labelled $i$. This exposed box is called {\it critical italic box}. 
In Fig. \sec.2, a circle is drawn around the label 
of each critical italic box.
On the other hand,
the $U$-coding is defined to be the set 
$\tmax(H(\l))$, where $H(\l)$ is the $H$-set of~$\l$.
A box in the leftmost column whose hook length is an element of the $U$-coding
is called {\it critical roman box}. In Fig. \sec.2, a circle is drawn
around the hook length number of each critical roman box.
Let us write the labels of all the exposed boxes (the vector $L=(L_i)$) with
its region numbers (the vector $R=(R_i)$)
and the $H$-set of $\l$ (the vector $H=(H_i)=H(\l)$), 
read from bottom to top.
{$$%
\font\cerclefont=cmsy10 at 14pt%
\def\a#1{#1\kern-9pt\hbox{\cerclefont\char'015}}%
\def\b#1{#1\kern-12pt\hbox{\cerclefont\char'015}}%
\def\c#1{#1\kern-11pt\hbox{\cerclefont\char'015}}%
\def\ng#1{\!\hbox{\rm -#1}}%
\matrix{
L=&\a{\it3} &\it 3 &\it 3 &\a{\it2}   &\a{\it4}&\it3&\it2&\it4&\it3&
\it 2& \it4& \it3& \it2& \a{\it1}&\a{\it0} & \cr
R=&\a{3}&2&1&\a{1}&\a0&0&0&\ng1&\ng1&\ng1&\ng2&\ng2&\ng2&\c{\ng2}&\c{\ng2}& \cr
H=&\b{23}&18&13&\b{12}&\a9&8&7&4&3&2&\ng1&\ng2&\ng3&\c{\ng4}&\c{\ng5}& \cr
}%
$$}%
It is easy to see that $L_j= (H_j-\ell(\l))\mod t$ 
and
$R_j= \lfloor(H_j-\ell(\l))/t\rfloor +1$.
This means that $L_i$ has a circle symbol if and only if $H_i$ has a 
circle symbol.
We then have a natural bijection
$$f: u_i \mapsto  \lfloor(u_i-\ell(\l))/t\rfloor +1 = n_{(u_i-\ell) \mod t}
\leqno{(\sec.7)}
$$
between the set $\{u_0, \ldots, u_{t-1}\}$ and 
$\{n_0, \ldots, n_{t-1}\}$. 
By (\sec.6) and (\sec.7) we have
$$
\leqalignno{
\{v_i'\} 
& = \{ tn_i+i-t'\} \cr
& = \{ t n_{ (u_i-\ell) \mod t} + (u_i-\ell) \mod t - t'\} \cr
& = \{ t( \lfloor(u_i-\ell)/t\rfloor +1  ) + (u_i-\ell) \mod t - t'\} \cr
& = \{ u_i-\ell+t'+1 \}. \cr
}
$$
On the other hand, $(v_i')$ is a $V$-coding,
because $v_i'\equiv i\mod t$ and 
$\sum v_i' = t\sum n_i + \sum i -t (t-1)/2 = 0$;
so that 
$$(\sum_i u_i )/t = \ell -t'-1. \leqno{(\sec.8)}$$
Hence
$$
\{v_i'\} 
 = \{ u_i-\ell+t'+1 \}
 = \{ u_i-(\sum_i u_i)/t \} 
 = \{ v_i \}.\qed 
$$

\medskip
Take again the same partition as in Example \sec.1; the $N$-coding
is
$$(n_0, n_1, n_2, n_3, n_4)=(-2, -2, 1, 3, 0).$$
We verify that
$$
\leqalignno{
&(v_0', v_1', v_2', v_3', v_4')\cr
&\qquad=(1\times5+0,\ 3\times5+1,\ 0\times5+2,\ -2\times5-2,\ -2\times5-1).\cr
&\qquad=(5, 16, 2, -12, -11)=(v_0, v_1, v_2, v_3, v_4).\cr
}
$$
\medskip

{\it Proof of (\sec.1) in Theorem \sec.1}.
From (\sec.6) we have
$$
\leqalignno{
\sum v_i^2&=\sum(tn_i+i-t')^2\cr
&=\sum \bigl((tn_i)^2 +2tin_i -2tt'n_i +i^2 +t'^2 -2it'  \bigr) \cr
&= t^2\sum n_i^2 + 2t\sum in_i + {(t-1)t(2t-1)\over 6} + tt'^2 - t' t (t-1)\cr
&= t^2\sum n_i^2 + 2t\sum in_i  + {t(t^2-1)\over 12}.
}
$$
Hence
$$
{1\over 2t}\sum v_i^2 = {t\over 2}\sum  n_i^2 + \sum in_i + {t^2-1\over 24} 
= |\l| + {t^2-1\over 24}.\qed
$$
\goodbreak\medskip\noindent
{\it \sec.4. Proof of (\sec.2)}. 
We first etablish the following two lemmas.
\proclaim Lemma \sec.5. 
For any $t$-compact set $A$ we have
$$
\prod_{a\in A, a>0} \Bigl(1-{t^2\over a^2} \Bigr)
=
\prod_{a\in \tmax(A), a\not=-t} {a+t\over a} .\leqno{(\sec.9)}
$$

\medskip

\noindent
{\it Example \sec.2}.
Take $t=5$. Then the set 
$$A=\{-5,-4, -3, -2, -1, 2,3,4,7,8,9,12,13,18,23\}$$ is $5$-compact.
We have $\tmax(A)=\{-5,-4, 9,12,23\}$. Hence
$$
\prod_{a\in A, a>0} \Bigl(1-{25\over a^2} \Bigr)=
{ 1\cdot 14 \cdot 17\cdot 28\over (-4) \cdot 9\cdot 12\cdot 23 }.
\leqno{(\sec.10)}
$$
{\it Proof}. Write
\vskip -10pt
$$
\prod_{a\in A, a>0} \Bigl(1-{t^2\over a^2} \Bigr)
=
\prod_{a\in A, a>0} {(a-t)\cdot (a+t)\over a \cdot a},
$$
then
delete the common factors in numerator and denominator, as illustrated 
by means of Example \sec.2. 
$$
{\def\coef#1|#2|#3|{{{#1 \over #3} \times {#2\over  #3}}}%
\def\ccf#1|#2|{{{#1 \over #2} }}%
\matrix{
\ccf1|-4| & & & & & & (a\equiv 1\mod 5)\cr
\noalign{\medskip}
\ccf2|-3|  & \coef-3|7|2| &\coef2|12|7| &\coef7|17|12| && &(a\equiv 2\mod 5)\cr
\noalign{\medskip}
\ccf3|-2| &\coef-2|8|3| &\coef3|13|8| &\coef8|18|13| &\coef13|23|18| 
   &\coef18|28|23| &(a\equiv 3\mod 5)\cr
\noalign{\medskip}
  \ccf4|-1|& \coef-1|9|4|  &\coef4|14|9| &&& &(a\equiv 4\mod 5)\cr
}}
$$
The product $(a-5)(a+5)/a^2$ for $a>0$ is reproduced in the row
determined by $a\mod 5$ in the above table, except for the leftmost column.
But the product of the factors in the leftmost column is equal to $1$ because
$t$ is an odd integer;
so that
the left-hand side of (\sec.10) is the product of the factors in the above
table. After deleting the common factors, it remains the rightmost 
fraction in each row.\qed

\medskip

\proclaim Lemma \sec.6.
Let $\l$ be a $t$-core and $(u_0, u_1, \ldots, u_{t-1})$ be its $U$-coding
(defined in the body of the construction of $\phi_V$).
Let $\l'$ be the $t$-core obtained from $\l$ by erasing the leftmost column
of $\l$ and $(u'_0, u'_1, \ldots, u'_{t-1})$ be its $U$-coding.
Then
$$
\prod_{0\leq i<j\leq t-1} {u_i-u_j \over u_i'-u_j'} = 
\prod_{j=1}^{t-1} {u_j+t\over u_j}. 
$$

\medskip
\noindent
{\it Example \sec.3}. Take the 5-core $\l$ given in Example \sec.1. The
$U$-coding of $\l$ is $(u_0,u_1,u_2,u_3,u_4)=(-5, -4, 12, 23, 9)$.
We have 
$$\l'=(13,9,5,5,3,3,3,1,1,1).$$ 
The $U$-coding of $\l'$ is
$(u'_0,u'_1,u'_2,u'_3,u'_4)=(-5, 11, 22, 8, -1)$.
Now, consider the cyclic rearrangement
$$(u''_0,u''_1,u''_2,u''_3,u''_4)=(-1, -5, 11, 22, 8)$$
 of $(u'_0,u'_1,u'_2,u'_3,u'_4)$.
We have $\prod(u'_i-u'_j)=\prod(u''_i-u''_j)$ because $t$ is an odd integer.
Moreover $u''_i=u_i-1$ for all $1\leq i\leq 4$. 
Hence
$$
\eqalignno{
 \prod_{0\leq i<j\leq t-1} { u_i-u_j \over u''_i-u''_j}
&=
 \prod_{j=1}^{t-1} { u_0-u_j \over u''_0-u''_j}\cr
&= { (-5+4)(-5-12)(-5-23)(-5-9) \over (-1+5)(-1-11)(-1-22)(-1-8) }\cr
&= { (-4+5)(12+5)(23+5)(9+5) \over (-4)(12)(23)(9) }.\cr
}
$$

{\it Proof}.
We suppose that $\l$ contains $\delta$ parts equal to $1$.
Its $H$-set  $H(\l)$ 
(viewed as a vector in decreasing order if necessary) 
can be split into six segments
$H(\l)=A_1A_2A_3A_4A_5A_6$ defined by (see Fig. \sec.3)

(i) $a \geq \delta+2$ for each $a\in A_1$;

(ii) $A_2=(\delta, \delta-1, \ldots, 3,2,1)$;

(iii) $A_3=(-1, -2, -3, \ldots, \delta+2-t)$;

(iv) $A_4=(\delta+1-t)$;

(v) $A_5=(\delta-t, \delta-1-t, \ldots, 1-t)$;

(vi) $A_6=(-t)$.

\smallskip

On the other hand the $H$-set $H(\l')$ of $\l'$
is split into five segments
$H(\l')=A_1'A_2'A_3'A_4'A_5'$ defined by 

(i') $A_1'=\{a-\delta-1 \ : \ a\in A_1\}$;

(ii') $A_2'=\{a-\delta-1 \ : \ a\in A_2\}=(-1, -2, \ldots, -\delta)$;

(iii') $A_3'=(-\delta-1)$;

(iv') $A_4'=\{a-\delta-1 \ : \ a\in A_3\}=(-\delta-2, -\delta-3, \ldots, -t+1)$;

(v) $A_5'=(-t)$.

%


\long\def\maplebegin#1\mapleend{}

\maplebegin

# --------------- begin maple ----------------------

# Copy the following text  to "makefig.mpl"
# then in maple > read("makefig.mpl");
# it will create a file "z_fig_by_maple.tex"

Hu:= 4;  # height unity
Hn:= 17; # height quantities
Lu:= 20; # large unity

X0:=-34; Y0:=0; # origin position

File:=fopen("z_fig_by_maple.tex", WRITE);

hline:=proc(x,y,len)
	fprintf(File, "\\hline(
end;

vline:=proc(x,y,len)
	fprintf(File, "\\vline(
end;

hdash:=proc(x,y,len)
local i, d,sp, xx;
	d:=1;
	sp:=1;
	for xx from x*Lu+X0 to x*Lu+X0+Lu*len-d by d+sp do
	fprintf(File, "\\hline(
	od:
end;

vdash:=proc(x,y,len)
local i, d,sp, yy;
	d:=1;
	sp:=1;
	for yy from y*Hu+Y0 to y*Hu+Y0+Hu*len-d by d+sp do
		fprintf(File, "\\vline(
	od:
end;

text:=proc(x,y,t)
	fprintf(File, "\\centerput(
			y*Hu+Y0+1-Hu,    t);
end;

textup:=proc(x,y,t)
	fprintf(File, "\\centerput(
			y*Hu+Y0+1-ceil(Hu/2),    t);
end;

# the partition

# h-line

hdash(0,18,1);
hdash(0,17,1); 
               hdash(1,14,1);
hdash(0,13,1); hdash(1,13,1);
hdash(0,12,1);
               hdash(1,9,1);
hline(0,8,1);  hdash(1,8,1);
hdash(0,4,1);  hline(1,4,1);
hdash(2,0,2);  hline(2,2,2);
hline(0,0,4);

# v-line

vdash(0,8,10);
                 vdash(1,8,10);
                 vline(1,4,4);  vdash(2,4,10);
                                vline(2,2,2);
vline(0,0,8);    vdash(1,0,4);  vdash(2,0,2);      vline(4,0,2);

# text
text(1,18, "$-t$");
text(1,17, "$1-t$");
text(1,16, "$\\vdots$");
text(1,15, "$\\delta-1-t$");
text(1,14, "$\\delta-t$");    text(2,14, "$-t$");
text(1,13, "$\\delta+1-t$");  text(2,13, "$1-t$");
text(1,12, "$\\delta+2-t$");  text(2,12, "$\\vdots$");
text(1,11, "$\\vdots$");      text(2,11, "$-\\delta-3$");
text(1,10, "$-2$");           text(2,10, "$-\\delta-2$");
text(1,9, "$-1$");            text(2,9, "$-\\delta-1$");
text(1,8, "$1$");             text(2,8, "$-\\delta$");
text(1,7, "$2$");             text(2,7, "$\\vdots$");
text(1,6, "$\\vdots$");       text(2,6, "$-2$");
text(1,5, "$\\delta$");       text(2,5, "$-1$");                  
#
text(0, 18, "$A_6$");
textup(0, 15, "$A_5$");  text(3, 14, "$A_5'$");
text(0, 13, "$A_4$");    textup(3, 11, "$A_4'$");
textup(0, 10, "$A_3$");  text(3, 9, "$A_3'$");
textup(0, 6, "$A_2$");   textup(3, 6, "$A_2'$");
textup(1, 2, "$A_1$");   textup(2, 2, "$A_1'$"); 
#                    

fclose(File);

# -------------------- end maple -------------------------
\mapleend


\newbox\boxarbre
\def\vline(#1,#2)#3|{\leftput(#1,#2){\lline(0,1){#3}}}
\def\hline(#1,#2)#3|{\leftput(#1,#2){\lline(1,0){#3}}}
\setbox\boxarbre=\vbox{\vskip
70mm\offinterlineskip 
%
\hline(-34,72)1|
\hline(-32,72)1|
\hline(-30,72)1|
\hline(-28,72)1|
\hline(-26,72)1|
\hline(-24,72)1|
\hline(-22,72)1|
\hline(-20,72)1|
\hline(-18,72)1|
\hline(-16,72)1|
\hline(-34,68)1|
\hline(-32,68)1|
\hline(-30,68)1|
\hline(-28,68)1|
\hline(-26,68)1|
\hline(-24,68)1|
\hline(-22,68)1|
\hline(-20,68)1|
\hline(-18,68)1|
\hline(-16,68)1|
\hline(-14,56)1|
\hline(-12,56)1|
\hline(-10,56)1|
\hline(-8,56)1|
\hline(-6,56)1|
\hline(-4,56)1|
\hline(-2,56)1|
\hline(0,56)1|
\hline(2,56)1|
\hline(4,56)1|
\hline(-34,52)1|
\hline(-32,52)1|
\hline(-30,52)1|
\hline(-28,52)1|
\hline(-26,52)1|
\hline(-24,52)1|
\hline(-22,52)1|
\hline(-20,52)1|
\hline(-18,52)1|
\hline(-16,52)1|
\hline(-14,52)1|
\hline(-12,52)1|
\hline(-10,52)1|
\hline(-8,52)1|
\hline(-6,52)1|
\hline(-4,52)1|
\hline(-2,52)1|
\hline(0,52)1|
\hline(2,52)1|
\hline(4,52)1|
\hline(-34,48)1|
\hline(-32,48)1|
\hline(-30,48)1|
\hline(-28,48)1|
\hline(-26,48)1|
\hline(-24,48)1|
\hline(-22,48)1|
\hline(-20,48)1|
\hline(-18,48)1|
\hline(-16,48)1|
\hline(-14,36)1|
\hline(-12,36)1|
\hline(-10,36)1|
\hline(-8,36)1|
\hline(-6,36)1|
\hline(-4,36)1|
\hline(-2,36)1|
\hline(0,36)1|
\hline(2,36)1|
\hline(4,36)1|
\hline(-34,32)20|
\hline(-14,32)1|
\hline(-12,32)1|
\hline(-10,32)1|
\hline(-8,32)1|
\hline(-6,32)1|
\hline(-4,32)1|
\hline(-2,32)1|
\hline(0,32)1|
\hline(2,32)1|
\hline(4,32)1|
\hline(-34,16)1|
\hline(-32,16)1|
\hline(-30,16)1|
\hline(-28,16)1|
\hline(-26,16)1|
\hline(-24,16)1|
\hline(-22,16)1|
\hline(-20,16)1|
\hline(-18,16)1|
\hline(-16,16)1|
\hline(-14,16)20|
\hline(6,0)1|
\hline(8,0)1|
\hline(10,0)1|
\hline(12,0)1|
\hline(14,0)1|
\hline(16,0)1|
\hline(18,0)1|
\hline(20,0)1|
\hline(22,0)1|
\hline(24,0)1|
\hline(26,0)1|
\hline(28,0)1|
\hline(30,0)1|
\hline(32,0)1|
\hline(34,0)1|
\hline(36,0)1|
\hline(38,0)1|
\hline(40,0)1|
\hline(42,0)1|
\hline(44,0)1|
\hline(6,8)40|
\hline(-34,0)80|
\vline(-34,32)1|
\vline(-34,34)1|
\vline(-34,36)1|
\vline(-34,38)1|
\vline(-34,40)1|
\vline(-34,42)1|
\vline(-34,44)1|
\vline(-34,46)1|
\vline(-34,48)1|
\vline(-34,50)1|
\vline(-34,52)1|
\vline(-34,54)1|
\vline(-34,56)1|
\vline(-34,58)1|
\vline(-34,60)1|
\vline(-34,62)1|
\vline(-34,64)1|
\vline(-34,66)1|
\vline(-34,68)1|
\vline(-34,70)1|
\vline(-14,32)1|
\vline(-14,34)1|
\vline(-14,36)1|
\vline(-14,38)1|
\vline(-14,40)1|
\vline(-14,42)1|
\vline(-14,44)1|
\vline(-14,46)1|
\vline(-14,48)1|
\vline(-14,50)1|
\vline(-14,52)1|
\vline(-14,54)1|
\vline(-14,56)1|
\vline(-14,58)1|
\vline(-14,60)1|
\vline(-14,62)1|
\vline(-14,64)1|
\vline(-14,66)1|
\vline(-14,68)1|
\vline(-14,70)1|
\vline(-14,16)16|
\vline(6,16)1|
\vline(6,18)1|
\vline(6,20)1|
\vline(6,22)1|
\vline(6,24)1|
\vline(6,26)1|
\vline(6,28)1|
\vline(6,30)1|
\vline(6,32)1|
\vline(6,34)1|
\vline(6,36)1|
\vline(6,38)1|
\vline(6,40)1|
\vline(6,42)1|
\vline(6,44)1|
\vline(6,46)1|
\vline(6,48)1|
\vline(6,50)1|
\vline(6,52)1|
\vline(6,54)1|
\vline(6,8)8|
\vline(-34,0)32|
\vline(-14,0)1|
\vline(-14,2)1|
\vline(-14,4)1|
\vline(-14,6)1|
\vline(-14,8)1|
\vline(-14,10)1|
\vline(-14,12)1|
\vline(-14,14)1|
\vline(6,0)1|
\vline(6,2)1|
\vline(6,4)1|
\vline(6,6)1|
\vline(46,0)8|
\centerput(-24,69){$-t$}
\centerput(-24,65){$1-t$}
\centerput(-24,61){$\vdots$}
\centerput(-24,57){$\delta-1-t$}
\centerput(-24,53){$\delta-t$}
\centerput(-4,53){$-t$}
\centerput(-24,49){$\delta+1-t$}
\centerput(-4,49){$1-t$}
\centerput(-24,45){$\delta+2-t$}
\centerput(-4,45){$\vdots$}
\centerput(-24,41){$\vdots$}
\centerput(-4,41){$-\delta-3$}
\centerput(-24,37){$-2$}
\centerput(-4,37){$-\delta-2$}
\centerput(-24,33){$-1$}
\centerput(-4,33){$-\delta-1$}
\centerput(-24,29){$1$}
\centerput(-4,29){$-\delta$}
\centerput(-24,25){$2$}
\centerput(-4,25){$\vdots$}
\centerput(-24,21){$\vdots$}
\centerput(-4,21){$-2$}
\centerput(-24,17){$\delta$}
\centerput(-4,17){$-1$}
\centerput(-44,69){$A_6$}
\centerput(-44,59){$A_5$}
\centerput(16,53){$A_5'$}
\centerput(-44,49){$A_4$}
\centerput(16,43){$A_4'$}
\centerput(-44,39){$A_3$}
\centerput(16,33){$A_3'$}
\centerput(-44,23){$A_2$}
\centerput(16,23){$A_2'$}
\centerput(-24,7){$A_1$}
\centerput(-4,7){$A_1'$}
%
%
}
\midinsert
$$
\box\boxarbre
$$
\centerline{\qquad\quad Fig. \sec.3. 
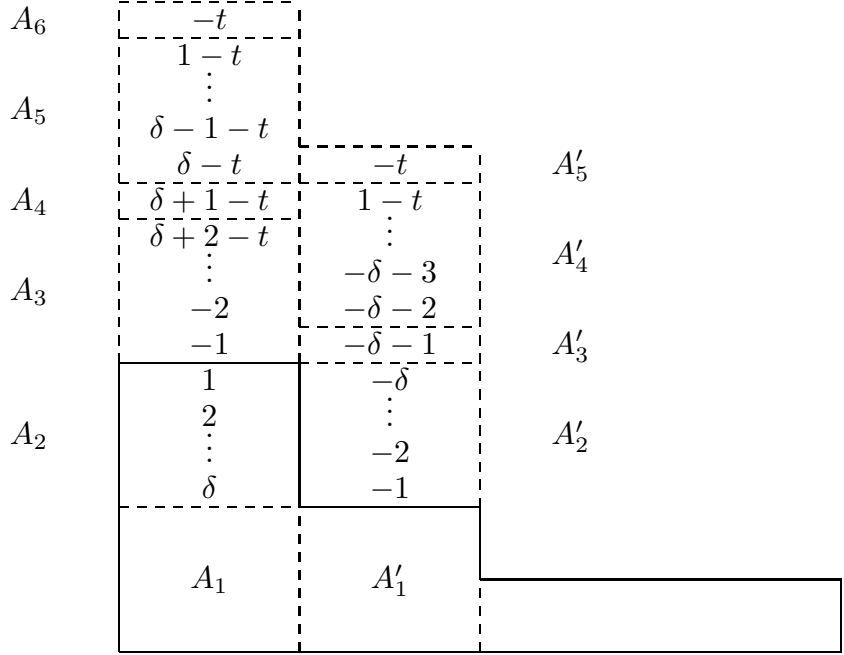
Comparison the hook lengths of $\l$ and $\l'$}
\endinsert
\smallskip
\noindent
Notice that some segments $A_i$ and $A_i'$ may be empty. More
precisely, 
$$
\cases{
A_2=A_5=A_2'=\emptyset, &  if $\delta=0$; \cr
A_3=A_4'=\emptyset, &  if $\delta=t-2$; \cr
A_3=A_4=A_3'=A_4'=\emptyset, & if $\delta=t-1$. \cr
}
$$

The basic facts are:

(i) $a\not \in \tmax(H(\l))$ for every $a\in A_5$; 
because 
$\{a\mod t\ : \ a\in A_5\} =\{a\mod t\ : \ a\in A_3\}$. 
In other words the set $A_5$  is {\it masked} by $A_3$. 

(ii) $\delta+1-t\in\tmax(H(\l))$; because $a\not\equiv 0\mod t$ for every 
$a\in A_1'$ so that
$a\not\equiv \delta+1 \mod t$ for every $a\in A_1$.
It is easy to see that
$a\not\equiv \delta+1 \mod t$ for every $a\in A_2\cup A_3$.

(iii) $-\delta-1\in\tmax(H(\l'))$; because $a\not\equiv 0\mod t$ for every
$a\in A_1 \cup A_2$ so that
$a\not\equiv -\delta-1\mod t$ for every $a\in A_1'\cup A_2'$.

(iv) Since that $a\mapsto a-\delta-1$ is a bijection between
$A_1\cup A_2\cup A_3$ and $A'_1\cup A'_2\cup A'_4$, 
it is also a bijection between
$\tmax(H(\l))\setminus\{-t,\delta-t+1\}$
and $\tmax(H(\l'))\setminus\{-t,-\delta-1\}$.
\medskip

The above facts enable us to derive the $U$-coding of $\l'$ from the
$U$-coding of $\l$ as follows. Let
$$(u_i)=(u_0=-t, u_1, u_2, \ldots, u_{k-1},  \delta+1-t, 
u_{k+1}, u_{k+1}, \ldots, u_{t-1})$$
be the $U$-coding of $\l$ and define
$$
(u''_i)=(u''_0=-\delta-1, u_1'', u_2'', \ldots, u_{k-1}'',  -t, 
u_{k+1}'', u_{k+1}'', \ldots, u_{t-1}'')$$
where $u''_i=u_i-\delta-1$ for $i\geq 1$.
Then, the $U$-coding of $\l'$ is simply
$$
(u'_i)=(u'_0=-t,  u_{k+1}'', u_{k+1}'', \ldots, u_{t-1}'',
-\delta-1, u_1'', u_2'', \ldots, u_{k-1}'').
$$
We have $\prod(u'_i-u'_j)=\prod(u''_i-u''_j)$ because $t$ is an odd integer.
On the other hand, 
$u''_i-u''_j=u_i-u_j$ for all $1\leq i<j\leq t-1$. Hence 
$$
\eqalignno{
 \prod_{0\leq i<j\leq t-1} { u_i-u_j \over u'_i-u'_j}
&= \prod_{0\leq i<j\leq t-1} { u_i-u_j \over u''_i-u''_j}
= \prod_{j=1}^{t-1} { u_0-u_j \over u''_0-u''_j}\cr
&= \prod_{j=1}^{t-1} { -t-u_j \over -\delta-1-u''_j}
= \prod_{j=1}^{t-1} { u_j +t\over u_j}.\qed\cr
}
$$
\medskip

{\it Proof of (\sec.2) in Theorem \sec.1}.
Because the $U$-coding and $V$-coding of~$\l$ only differ by the
normalization given in (\sec.3) and $t$ is an odd integer, we have
$\prod(v_i-v_j)=\prod(u_i-u_j)$.
By Lemmas \sec.6 and \sec.5 we have 
$$
\leqalignno{
\prod_{0\leq i<j\leq t-1} (u_i-u_j) 
&=\prod_{j=1}^{t-1} {u_j+t\over u_j}\times 
  \prod_{0\leq i<j\leq t-1} (u'_i-u'_j) \cr
&= 
\prod_{a\in H(\l), a>0} \Bigl(1-{t^2\over a^2} \Bigr)\times
  \prod_{0\leq i<j\leq t-1} (u'_i-u'_j) \cr
&= \cdots 
=
K\times \prod_{v\in \l} \Bigl(1-{t^2\over h_v^2} \Bigr). \cr
}
$$
Taking $\l$ as the empty $t$-core, the $U$-coding of $\l$ is
$(-t, -t+1, -t+2, \ldots, -3, -2, -1)$. 
We then obtain $K= (-1)^{t'}{ 1!\cdot 2!\cdot 3!\cdots (t-1)!} $\qed
\medskip \noindent
{\it \sec.5. Enf of the proof of the Main Theorem.} 
Recall that the Dedekind $\eta$-function is defined by
$$
\eta(x)=x^{1/24} \prod_{m\geq 1} (1-x^m). \leqno{(\sec.11)}
$$
Let $t=2t'+1$ be an odd integer. Macdonald obtained the following result 
[Ma72] (see comments in \S 1).

\proclaim Theorem \sec.7 [Macdonald].
We have
$$
\eta(x)^{t^2-1} = c_0 \sum_{(v_0, \ldots, v_{t-1})} \prod_{i<j} (v_i-v_j) 
x^{(v_0^2+v_1^2+\cdots+v_{t-1}^2)/(2t)},
\leqno{(\sec.12)}
$$
where the sum ranges over all $V$-codings $(v_0, v_1, \ldots, v_{t-1})$ 
(see Definition \sec.1)
and $c_0$ is a numerical constant.

Consider the term of lowest degree in the above power series.
We immediately get 
$$
c_0={(-1)^{t'}\over 1!\cdot 2!\cdot 3!\cdots (t-1)!}. \leqno{(\sec.13)} 
$$

{\it Proof of the Main Theorem}.
Let $n\geq 0$ be a positive integer.
The coefficient $C_n(\beta)$ of $x^n$ on the left-hand side of the 
Main Identity is a polynomial
in $\beta$ of degree $n$. The coefficient $D_n(\beta)$ of $x^n$ on the 
right-hand side
of the Main Identity is also a 
polynomial in $\beta$ of degree $n$ thanks to (2.3). 
For  proving $C_n(\beta)=D_n(\beta)$, it suffices to find $n+1$
explicit numerical values $\beta_0, \beta_1, \ldots, \beta_n$ such that 
$C_n(\beta_i)= D_n(\beta_i)$ for $0\leq i\leq n$
by using the Lagrange interpolation formula.
The basic fact is that
$$
\prod_{v\in\l}\bigl(1-{t^2 \over h_v^2}\bigr)=0
$$
for every partition $\l$ which is not a $t$-core.
By comparing Theorems \sec.1 and \sec.7 we see that the Main Identity is true
when $\beta=t^2$ for every odd integer~$t$, i.e.,
$$
\sum_{\l\in \setP}\ x^{|\l|} \prod_{v\in\l}\bigl(1-{t^2 \over h_v^2}\bigr)
\ =\ \prod_{m\geq 1} {(1-x^m)^{t^2-1}},
$$
so that $C_n(\beta)=D_n(\beta)$ for every complex number $\beta$.  \qed
\medskip

Note that Kostant already observed that $C_n(\beta)$ 
is a polynomial in~$\beta$, but did not mention any explicit expression [Ko04].



\def\sec{6}
\section{\sec. New formulas about hook lengths} 

In Sections 2-4 we have taken special numerical values for $\beta$.
In this section we  compare the coefficients of $\beta^k$
to derive further identities. 
\goodbreak\medskip \noindent
{\it \sec.1. Comparing the coefficients of $\beta$}. 
\proclaim Proposition \sec.1.
We have
$$
\sum_{n\geq 1}x^n\sum_{\l\vdash n}\sum_{v\in\l} 
{1\over h_v^2}
=\prod_{m\geq 1} {1\over 1-x^m} \times 
\sum_{k\geq 1} {x^k\over k (1-x^k)}. \leqno{(\sec.1)}
$$

{\it Proof}. 
Identity (\sec.1) follows from (2.3) by comparing the coefficients 
of~$\beta^1$ on both sides. \qed
\medskip \noindent
{\it \sec.2. The Stanley-Elder-Bessenrodt-Bacher-Manivel Theorem}. 
We also have a second proof of Proposition \sec.1 that is direct and 
provides a more general result about the power sum 
of the hook lengths.

\proclaim Theorem \sec.2.
We have
$$
\sum_{n\geq 1}x^n\sum_{\l\vdash n}\sum_{v\in\l} 
{ h_v^\alpha}
=\prod_{m\geq 1} {1\over 1-x^m} \times 
\sum_{k\geq 1} {x^k k^{\alpha+1}\over 1-x^k}. \leqno{(\sec.2)}
$$

Recall that $\sigma_\alpha(n)=\sum_{d|n} d^\alpha$ is the 
$\alpha$-th power sum of all positive divisors of $n$ (see [Se70, p.149])
whose generating function is classically given by
$$
\sum_{k\geq 1}{x^k k^\alpha\over 1-x^k} = \sum_{n\geq 1} \sigma_\alpha(n)x^n.
\leqno{(\sec.3)}
$$
Using (\sec.3)
identity (\sec.2) can be rewritten as
$$
\sum_{n\geq 1}x^n\sum_{\l\vdash n}\sum_{v\in\l} 
{ h_v^\alpha}
=\prod_{m\geq 1} {1\over 1-x^m} \times \sum_{n\geq 1} 
\sigma_{\alpha+1}(n)x^n.  \leqno{(\sec.4)}
$$
\medskip

The proof of Theorem \sec.2 is based on an elegant result about
the multi-set of hook lengths and the multi-set of parts of all partitions
of~$n$. 
Many studies have been done along those lines 
[Be98, BM02, Ho86, St04, KS82, We1, We2]. 
Each hook length $h_v$ can be split into $h_v=a_v+l_v+1$ where
$a_v$ is the {\it arm length} and $l_v$ is the {\it leg length} 
(see [St99, p.457]).
The ordered pair $(a_v, l_v)$ is called a {\it hook type}. 

\proclaim Theorem \sec.3 [Stanley-Elder-Bessenrodt-Bacher-Manivel].
\quad 
Let $n\geq k \geq 1$ be two integers. Then for every positive $j<k$
the total number of occurrences of the part $k$ among all partitions of $n$
is equal to the number of boxes whose hook type is $(j, k-j-1)$.

We now state a weaker form of the SEBBM Theorem, much easier to figure out.
Let $A$ be a multi-set of positive integers. Define $\dot A$ to be
the multi-set derived from $A$ by replacing each element $a$ of $A$
by $a$ copies of $a$.
For instance, with $A=\{1,1,2,5\}$
we obtain $\dot A=\{1,1,2,2,5,5,5,5,5\}=\{1^1,1^1,2^2,5^5\}$.

\proclaim Proposition \sec.4.
Let $H(n)$ (resp. $G(n)$) be the multi-set of 
all hook lenghts (resp. the parts) of all partitions of $n$. 
Then
$$
H(n)=\dot G(n).\leqno{(\sec.5)}
$$

For example, the set of all partitions of $4$ with their 
hook length multi-sets is reproduced in Fig. \sec.1.
We see that 
$H(4)=\{1^7,2^6,3^3,4^4\}$.
On the other hand, $G(4)=\{1,1,1,1,2,1,1,2,2,3,1,4\}$.
We have $\dot G(4)=\{1^7,2^6,3^3,4^4\}=H(4)$.
Notice that $\dot G(4)$ can be represented as in Fig.~\sec.2.

{
\medskip
\setbox1=\hbox{$\ytableau{12pt}{0.4pt}{}
{\\1   \cr 
 \\2   \cr
 \\3   \cr
 \\4   \cr
}$}
\setbox2=\hbox{$\ytableau{12pt}{0.4pt}{}
{\\1    \cr 
 \\2    \cr
 \\4 &\\1  \cr
}$}
\setbox3=\hbox{$\ytableau{12pt}{0.4pt}{}
{\\2 &\\1    \cr
 \\3 &\\2  \cr
}$}
\setbox4=\hbox{$\ytableau{12pt}{0.4pt}{}
{\\1    \cr 
 \\4 &\\2 &\\1  \cr
}$}
\setbox5=\hbox{$\ytableau{12pt}{0.4pt}{}
{ \\4 &\\3 &\\2 &\\1  \cr
}$}
\centerline{\box1\qquad\box2\qquad\box3\qquad\box4\qquad\box5}
\nobreak \smallskip \nobreak
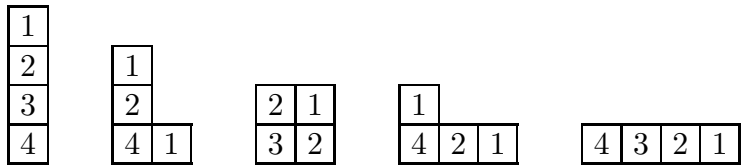
\centerline{Fig. \sec.1. The multi-set $H(4)$ of hook lengths}
\medskip
}

{
\setbox1=\hbox{$\ytableau{12pt}{0.4pt}{}
{\\1   \cr 
 \\1   \cr
 \\1   \cr
 \\1   \cr
}$}
\setbox2=\hbox{$\ytableau{12pt}{0.4pt}{}
{\\1    \cr 
 \\1    \cr
 \\2 &\\2  \cr
}$}
\setbox3=\hbox{$\ytableau{12pt}{0.4pt}{}
{ \\2 &\\2    \cr
 \\2 &\\2  \cr
}$}
\setbox4=\hbox{$\ytableau{12pt}{0.4pt}{}
{\\1    \cr 
 \\3 &\\3 &\\3  \cr
}$}
\setbox5=\hbox{$\ytableau{12pt}{0.4pt}{}
{ 
\\4 &\\4 &\\4 &\\4  \cr
}$}
\centerline{
\box1\qquad\box2\qquad\box3\qquad\box4\qquad\box5
}
\nobreak \smallskip \nobreak
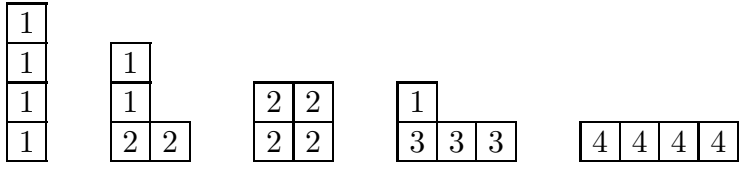
\centerline{Fig. \sec.2. The multi-set $\dot G(4)$ of parts with duplications}
\medskip
}

Theorem \sec.3 can be used to evaluate the power sum of the hook lengths.
We obtain immediately the following Theorem, which means that
the $r$-th power sum of the hook lengths is equal to the $(r+1)$-st
power sum of the parts.

\proclaim Corollary \sec.5.
For each positive integer $n$ and each  complex number $\alpha$  we have
$$
\sum_{\l\vdash n}\sum_{v\in\l} h_v^\alpha
=
\sum_{\l\vdash n}\sum_{i} \l_i^{\alpha+1}.\leqno{(\sec.6)}
$$
\medskip

By Corollary \sec.5, we see that Theorem \sec.2 is equivalent to the next 
Theorem.

\proclaim Theorem \sec.6.
We have
$$
\sum_{n\geq 1}x^n \sum_{\l\vdash n}\sum_{i} \l_i^{\alpha}
=\prod_{m\geq 1} {1\over 1-x^m} \times 
\sum_{k\geq 1} {x^k k^{\alpha}\over 1-x^k}. \leqno{(\sec.7)}
$$

{\it Proof}.
Let
$$
F(k):=\sum_{n_k, n_{k+1},  \ldots\geq 0} 
x^{kn_k+(k+1)n_{k+1} +\cdots}
(n_k k^\alpha+n_{k+1} (k+1)^\alpha + \cdots).
\leqno{(\sec.8)}
$$
We have
$$
\leqalignno{
F(k)
&= \sum_{n_k} x^{k n_k} n_k k^\alpha \times 
  \sum_{n_{k+1},\ldots} x^{(k+1)n_{k+1}+\cdots}
+ \sum_{n_k} x^{kn_k} F(k+1)\cr
&=
{k^\alpha \over (1-x^{k+1}) (1-x^{k+2})\cdots } \sum_n x^{kn} n +
{1\over 1-x^k} F(k+1)\cr
&=
{1 \over (1-x^{k}) (1-x^{k+1})\cdots } {k^\alpha x^k\over 1-x^k} +
{1\over 1-x^k} F(k+1). &(\sec.9)\cr
}
$$
Let
$$
F'(k)={1\over (1-x)(1-x^2)\cdots (1-x^{k-1})}F(k).\leqno{(\sec.10)}
$$
Then identity (\sec.9) becomes
$$
F'(k) =
{1 \over (1-x) (1-x^2)\cdots } {k^\alpha x^k\over 1-x^k} +
F'(k+1).\leqno{(\sec.11)}
$$
By iteration 
$$
F'(1)=\prod_{m\geq 1} {1\over 1-x^m} \times 
\sum_{k\geq 1} {x^k k^{\alpha}\over 1-x^k}. \leqno{(\sec.12)}
$$
Thus, the left-hand side of (\sec.7) is equal to $F(1)=F'(1)$. \qed

\medskip
Putting $\alpha=0$ and $\alpha=-1$ we obtain the following specializations.
\proclaim Corollary \sec.7.
We have
$$
\leqalignno{
\sum_{n\geq 1}x^n\sum_{\l\vdash n}\sum_{v\in\l} 
{ 1}
&=\prod_{m\geq 1} {1\over 1-x^m} \times 
\sum_{k\geq 1} {x^k k\over 1-x^k}\cr
&=x{d\over dx}\prod_{m\geq 1} {1\over 1-x^m}   &{(\sec.13)}\cr
}
$$
\vskip -10pt
\noindent
and
\vskip -10pt
$$
\leqalignno{
\sum_{n\geq 1}x^n\sum_{\l\vdash n}\sum_{v\in\l} 
{1\over h_v}
&=\prod_{m\geq 1} {1\over 1-x^m} \times 
\sum_{k\geq 1} {x^k \over 1-x^k}\cr
&= \sum_{m\geq 1} {mx^m\over (1-x)(1-x^2)\cdots (1-x^m)}. &(\sec.14)\cr
}
$$

\noindent
For the second equality of (\sec.14), see [Slo, A006128].

\medskip
{\it Historical Remarks about the SEBBM theorem}.
In the present paper we do not give the proof of the SEBBM theorem. 
We only want to 
make the following historical remarks. Stanley proved the case $j=0$ in 1972.
Independent discoveries and proofs were given by Kirdar and
Skyrme (1982), Paul Elder (1984) and Hoare (1986) 
(see [We1, We2, St04, KS82, Ho86]).
This result is called Elder's Theorem.
Bessenrodt [Be98] proved the general case of Theorem \sec.3 in 1998.
The final version of this result was given by Bacher and 
Manivel [BM02] in 2002.
In fact, when we re-discovered Theorem \sec.3, as will be further explained, 
we noticed that Elder's Theorem, stated in this hook length language,
was just the particular case $j=0$ of Theorem \sec.3.

When preparing the present paper we rediscovered Theorem \sec.3 in the 
following manner. First, we obtained Proposition \sec.1, as mentioned earlier 
by comparing the coefficients of $\beta$ in the Main Theorem.
We then expanded the right-hand side of (\sec.1) and calculated 
the first terms: 
$$
\sum_{n\geq 1}x^n\sum_{\l\vdash n}\sum_{v\in\l} 
{1\over h_v^2}
=
x + 5{ x^2\over 2!} + 29{x^3\over 3!} + 218 {x^4\over 4!} + 
  1814 {x^5\over 5!} +\cdots
\leqno{(\sec.15)}
$$
When searching for the sequence $1,5,29,218,1814,\ldots$ in {\it
The On-Line Encyclopedia of Integer Sequences} [Slo]  we got the sequence
A057623, that referred to
``{\tt $n!$ * (sum of reciprocals of all parts in 
unrestricted partitions of $n$).}"
Next we calculated
$$
\sum_{n\geq 1}x^n\sum_{\l\vdash n}\sum_{v\in\l} 
{1\over h_v}=x+3x^2+6x^3+12x^4+20x^5+35x^6+\cdots
$$
by enumerating all partitions. Going back to the
The On-Line Encyclopedia the sequence $(1,3,6,12,20,35,\ldots)$ referred to
the sequence A006128 with the following
description:
``{\tt Total number of parts in all partitions of $n$.}"
Those facts led us to discover equality (\sec.6), and then Theorem \sec.3.

\medskip \noindent
{\it \sec.3. Comparing the coefficients of $\beta^2$}. 
By selecting the coefficients of $\beta^2$ in our Main Identity
we obtain the following
equality about hook lengths. Unlike Theorem \sec.2 the following
results can not be derived from the SEBBM theorem.

\proclaim Proposition \sec.8.
We have
$$
\sum_{n\geq 2}x^n\sum_{\l\vdash n}\sum_{\{u,v\}} 
{1\over h_u^2 h_v^2}
={1\over 2}\prod_{m\geq 1} {1\over 1-x^m} \times 
\Bigl(\sum_{k\geq 1} {x^k\over k (1-x^k)}\Bigr)^2, \leqno{(\sec.16)}
$$
where the third sum ranges over all unordered pairs $\{u,v\}$ such that 
$u,v\in\l$ and $u\not= v$. 

By Theorem \sec.7 and Proposition \sec.1, we have
$$
\sum_{n\geq 2}x^n\sum_{\l\vdash n}\sum_{(u,v)} 
{1\over h_u^2 h_v^2}
=
\prod_{m\geq 1} (1-x^m) \times 
\Bigl(\sum_{n\geq 1}x^n\sum_{\l\vdash n}\sum_{v\in\l} {1\over h_v^2}\Bigr)^2,
\leqno{(\sec.17)}
$$
where the third sum ranges over all ordered pairs $(u,v)$ such that 
$u,v\in\l$ and $u\not= v$. 

\proclaim Theorem \sec.9 [=1.4].
We have
$$
\sum_{n\geq 1} x^n \sum_{\l\vdash n} 
\bigl(\sum_{v\in\l} {1\over h_v^2}\bigr)^2
=
\prod_{m\geq 1}{1\over 1-x^m} \Bigl(
\sum_{k\geq 1}{x^k k^{-3}\over 1-x^k}+
\bigl(\sum_{k\geq 1}{x^k k^{-1}\over 1-x^k}\bigr)^2
\Bigr). \leqno{(\sec.18)}
$$

\medskip
{\it Proof}. For each partition $\l$ we have
$$
\bigl(\sum_{v\in\l} {1\over h_v^2}\bigr)^2
=
\sum_{v\in\l} {1\over h_v^4}
+
2\sum_{\{u, v\}} {1\over h_u^2 h_v^2}
$$
and conclude in view of Theorem \sec.2 and Proposition \sec.8.\qed

\medskip \noindent
{\it \sec.4. Comparing the coefficients of $\beta^{n}x^n$  
             and $\beta^{n-1}x^n$}.  
Recall that $f_\l$ is the number of standard Young tableaux of shape $\l$.
Comparing the coefficients
of $(-\beta)^n x^n$ on both sides of the Main Theorem (see, for example (2.3)),
we get
$$
\sum_{\l\vdash n} f_\l^2 =  n! \leqno{(\sec.19)}
$$

\proclaim Theorem \sec.10 [=1.3, marked hook formula].
We have
$$
\sum_{\l\vdash n} f_\l^2\sum_{v\in\l} h_v^2 = {n(3n-1)\over 2} n!
\leqno{(\sec.20)}
$$

{\it Proof}.
Selecting the coefficients of $(-\beta)^{n-1}x^n$ on the right-hand side of 
equation~(2.3) we obtain
$$
\leqalignno{
&[(-\beta)^{n-1}x^n]\ 
\prod_{m\geq 1} {1\over 1-x^m}\times 
\exp{\Bigl(-\beta\sum_{k\geq 1} {x^k\over k(1-x^k)}\Bigr)} \cr
=
&[x^n]\ 
{1\over (n-1)!} \prod_{m\geq 1} {1\over 1-x^m}\times 
\Bigl(\sum_{k\geq 1} {x^k\over k(1-x^k)}\Bigr)^{n-1} \cr
=
&[x^1]\  
{1\over (n-1)!}  {1\over 1-x}\times 
\Bigl( {1\over 1-x}+{x\over 2(1-x^2)}\Bigr)^{n-1} \cr
=& {n(3n-1)\over 2 n!}. &{(\sec.21)}\cr
}
$$
Selecting the coefficients of $(-\beta)^{n-1}x^n$ on the left-hand side of 
equation~(2.3) we get
$$
\sum_{\l\vdash n}\ \sum_{u\in\l} \prod_{v\not=u}{1 \over h_v^2}
=
\sum_{\l\vdash n}\ \prod_{u\in\l}{1 \over h_u^2} \sum_{v\in\l}  h_v^2
=
\sum_{\l\vdash n}\ {f_\l^2\over n!^2} \sum_{v\in\l}  h_v^2. \qed
\leqno{(\sec.22)}
$$

{\it Remark}.
Is there a combinatorial proof of the marked hook formula (\sec.20),
analogous to the Robinson-Schensted-Knuth correspondence 
for proving (\sec.19)~?
Let $T$ be a standard Young tableau of shape $\l$ (see [Kn98, p.47]),
$u$ be a box in $\l$ and 
$m$ an integer such that $1\leq m\leq h_u(\l)$. 
The triplet $(T,u,m)$ is called
a {\it marked Young tableau} of shape $(\l, u)$.
The
number of marked Young tableaux of shape $(\l, u)$ is then $f_\l h_u$.
On the other hand,
call {\it marked permutation} each triplet $(\sigma, j,k)$ where 
$\sigma\in\setS_n$, $1\leq j\leq n$ and $1\leq k\leq n+j-1$. We say that 
the letter $j$ within the permutation $\sigma$ is marked $k$.
The total number of marked permutations of order $n$ is
$$
\sum_{j=1}^n (n+j-1) n!= {n(3n-1)\over 2} n!
$$
\noindent
{\it Example}. The sequence $6\;4\;9\;5_k\;7\;1\;2\;8\;3$ 
with $1\leq k\leq 13$ is a marked permutation, 
whose letter $5$ is marked $k$.
The following two diagrams
are two marked Young tableaux of the same shape, where $1\leq i,j\leq 3$.
\vskip -5pt
{
\setbox1=\hbox{$
\def\b{\\{\hbox{}}}
\ytableau{14pt}{0.4pt}{}
{\\5 &\\9  \cr 
 \\2 &\\8  \cr
 \\1 &\\4 &\\7   \cr
\noalign{\vskip 3pt}
}
\centerput(-11,7.3){${}_i$}
$
}
\setbox2=\hbox{$
\ytableau{14pt}{0.4pt}{}
{\\8 &\\9  \cr 
 \\4 &\\5  \cr
 \\1 &\\3 &\\7   \cr
\noalign{\vskip 3pt}
}
\centerput(-11,7.7){${}_j$}
$
}
$$\box1\qquad\qquad\box2$$
}
\vskip -10pt
\noindent
For proving the marked hook formula we need 
find a {\it marked} Robinson-Schensted-Knuth correspondence between
pairs of marked Young tableaux and marked permutations. 

\medskip \noindent
{\it \sec.5. Comparing the coefficients of $\beta^{n-2}x^n$   
             and $\beta^{n-3}x^n$}.
In the same manner as in the proof of the marked hook formula 
we obtain the following results
by selecting the coefficients of $(-\beta)^{n-2}x^n$
and  $(-\beta)^{n-3}x^n$ in (2.3). 
\proclaim Proposition \sec.11.
We have
$$
\sum_{\l\vdash n} f_\l^2\sum_{\{u,v\}} h_u^2 h_v^2 = 
{n(n-1)(27n^2-67n+74)\over 24} n!,
$$
where the second sum ranges over all unordered pairs $\{u,v\}$ such that 
$u,v\in\l$ and $u\not= v$. 

\proclaim Proposition \sec.12.
We have
$$
\sum_{\l\vdash n} f_\l^2\sum_{\{u,v,w\}} h_u^2 h_v^2 h_w^2 = 
{n(n-1)(n-2)(27n^3-174n^2+511n-600)\over 48} n!,
$$
where the second sum ranges over all unordered distinct triplets 
$\{u,v,w\}$ of boxes of the partition $\l$.

\def\sec{7}
\section{\sec. Improvement of a result due to Kostant} 
Let 
$$
\prod_{n\geq 1} (1-x^n)^s  =  \sum_{k\geq 0} f_k(s) x^k. \leqno{(\sec.1})
$$
Kostant proved the following result [Ko04, Th. 4.28]. 

\proclaim Theorem \sec.1 [Kostant].
Let $k$ and $m$ be two positive integers such that $m\geq \max(k,4)$. 
Then $f_k(m^2-1)\not=0$. 

The condition $m>1$ in the original Theorem of Kostant should be replaced by
$m\geq 4$, as, for example, we have $f_3(8)=0$ (see Theorem \sec.3).

\proclaim Theorem \sec.2 [=1.5].
Let $k$ be a positive integer and $s$ be a real number 
such that $s\geq k^2-1$.  Then $(-1)^k f_k(s)>0$. 

\smallskip
\noindent
{\it Remarks}. We extend Kostant's result in two directions: first,
we claim that $(-1)^k f_k(s)>0$ instead of $f_k(s)\not=0$; second,
$s$ is any real number instead of an integer of the form $m^2-1$. 

\smallskip
{\it Proof}.
By the Main Theorem we may write 
$$
(-1)^kf_k(s)=\sum_{\l\vdash k} W(\l), 
\leqno{(\sec.2)}
$$
where
$$
W(\l)= \prod_{v\in\l} \Bigl({s+1\over h_v^2}-1\Bigr)
= \prod_{v\in\l} \Bigl({s+1-h_v^2\over h_v^2}\Bigr). 
\leqno{(\sec.3)}
$$
For each $\l\vdash k$ and $v\in\l$ we have $h_v(\l)\leq k$, so that
$W(\l)\geq 0$.  This means that there is {\it no cancellation} in 
the sum  (\sec.2).
If $s>k^2-1$ then $W(\l)>0$. If  $s=k^2-1\geq 15$ we have $k\geq 4$.
In that case there is at least one partition $\l$, whose hook lengths
are strictly less than $k$. Hence $W(\l)>0$. \qed 
\medskip

Here is another result of Kostant [Ko04, Th.4.27].
\proclaim Theorem \sec.3 [Kostant].
We have
$$
\leqalignno{
f_4(s)&=1/4!\ s(s-1)(s-3)(s-14); \cr
-f_3(s)&=1/3!\ s(s-1)(s-8);\cr
f_2(s)&=1/2!\ s(s-3).\cr
}
$$

Even though we do not see how to factorize each $f_k(s)$, the occurrences 
of some factors in the above formulas have some relevance in terms of 
hook lengths.
Every partition contains one hook length $h_v=1$, so that $f_k(s)$ has
the factor $s+1-h_v^2=s$ (see (\sec.3)). Every partition of $3$ contains
a hook length $h_v=3$, so that $f_3(s)$ has the factor $s-8$. Every
partition of $2$ or $4$ has a hook length $h_v=2$, so that 
$s-3$ is a factor of $f_2(s)$ and $f_4(s)$.


\def\sec{8}
\section{\sec. The magic partition formula} 
Let $s_\l$ be the Schur functions corresponding to the partition $\l$ (see 
[Ma95, p.40; St99, p.308; La01, p.8]). 
Let $X=\{x_1, x_2, \ldots\}$ and $Y=\{y_1, y_2, \ldots\}$ be two alphabets.
The famous Cauchy formula is stated as follows 
(see [Ma95, p.63; St99, p.322; La01, p.13; Kn70]):
\proclaim Theorem \sec.1 [Cauchy].
$$
\prod_{i,j} {1\over 1-x_i y_j} = \sum_{\l\in\setP} s_\l(X) s_\l(Y).
\leqno{(\sec.1)}
$$

Let $d$ be a positive integer. Taking 
$$
X=\{x, x^2, x^3, \ldots\} \hbox{\quad and\quad} 
Y=\{1,1,\ldots, 1\} = \{1^d\},
$$
we get the following specialization:
$$
\prod_{m\geq 1} \Bigl({1\over 1-x^m}\Bigr)^d 
= \sum_{\l\in\setP} s_\l(x,x^2,\ldots) s_\l(1^d).
\leqno{(\sec.2)}
$$
Also recall the classical hook-content formula [Ro58; St99, p.374]:
\proclaim Theorem \sec.2.
For any partition $\l$ and positive integer $n$ we have
$$
s_\l(1,x,x^2,\ldots, x^{n-1})=x^{b(\l)} 
\prod_{v\in\l} {1-x^{n+c(v)}\over 1-x^{h_v}}
\leqno{(\sec.3)}
$$
where $b(\l)=\sum_i (i-1) \l_i$ and $c(v)=j-i$ if $v=(i,j)\in\l$.

By Theorem \sec.2 we then have 
$$
\leqalignno{
s_\l(x,x^2,x^3, \ldots)
&=x^{|\l|+b(\l)} \prod_{v\in\l} {1\over 1-x^{h_v}}, &{(\sec.4)}\cr
s_\l(1^d)
&= \prod_{v\in\l} {d+c(v) \over h_v}. &{(\sec.5)}\cr
}
$$

\proclaim Theorem \sec.3.
For any complex number $\beta$ we have
$$
\prod_{m\geq 1} (1-x^m)^\beta
=
\sum_{\l\in\setP} x^{|\l| +b(\l)} 
\prod_{v\in\l} { c(v)-\beta\over h_v (1-x^{h_v}) }. \leqno{(\sec.6)}
$$

{\it Proof.}
From (\sec.2), (\sec.4) and (\sec.5) we see that  
(\sec.6) is true for any negative integer $\beta$.
Thus (\sec.6) is true for any complex number $\beta$ 
(see the explanation given in the proof of the Main Theorem, \S5.5). \qed

\medskip

{\it Remark \sec.1}.
Theorem \sec.3 appears to be another formula for all the powers of 
the Euler Product. Although its form is 
analogous with the Main Theorem, 
it has fewer applications. As the 
variable $x$ occurs in the denominator
on the right-hand side of (\sec.6), 
it becomes cumbersome to select specific coefficients of $x^n$. 
There are apparently no specialization leading to
Macdonald identities.
\medskip

{\it Remark \sec.2}.
When $\beta$ is given the value 1 in (\sec.6) 
we recover the following identity due to 
Euler [An76, p.11].

\proclaim Corollary \sec.4 [Euler].
We have
$$
\prod_{m\geq 1} (1-x^m) = 
\sum_{n\geq 0} {(-1)^n x^{n(n+1)/2}\over (1-x)(1-x^2)\cdots (1-x^n)}.
$$

Combining Thorem \sec.3 and the Main Theorem we get the following
result, called {\it magic partition formula} because the sum and the product on 
both sides range over the same sets $\l\in\setP$ and $v\in\l$. 

\proclaim Theorem \sec.5 [Magic partition formula].
For any complex number~$\beta$ we have
$$
\sum_{\l\in\setP} x^{|\l| +b(\l)} 
\prod_{v\in\l} { c(v)+1-\beta \over h_v (1-x^{h_v}) }
\ =\ 
\sum_{\l\in \setP}\ x^{|\l|}\prod_{v\in\l}{h_v^2-\beta \over h_v^2}.
\leqno{(\sec.7)}
$$


\def\sec{9}
\section{\sec. Reversion of the Euler Product} 
Let $y(x)$ be a formal power series satisfying the following relation
$$
\leqalignno{
x&=y(1-y)(1-y^2)(1-y^3)\cdots &{(\sec.1)}\cr
&=y-y^2-y^3+y^6+y^8-y^{13}-y^{16}+\cdots\cr
}
$$
The first coeficients of the reversion series in (\sec.1)
are the following
$$
y(x)=
x+{x}^{2}+3\,{x}^{3}+10\,{x}^{4}+38\,{x}^{5}+153\,{x}^{6}+646\,{x}^{7}
+\cdots \leqno{(\sec.2)} 
$$
They are referred to as the first values of the sequence A109085
in The On-Line Encyclopedia of Integer Sequences [Slo].

\proclaim Theorem \sec.1.
We have the following explicit formula for the reversion of (\sec.1)
in terms of hook lengths:
$$
y(x)=\sum_{n\geq 1}{x^n\over n} \sum_{\l\vdash n-1}\ \prod_{v\in\l}
\bigl(1+{n-1\over h_v^2} \bigr). \leqno{(\sec.3)}
$$

{\it Proof}. Rewrite (\sec.1) as $y=x\phi(y)$ where 
$\phi(y)=\prod_{m\geq 1} (1-y^m)^{-1}$.
By the Lagrange inversion formula and the Main Theorem we have
$$
\leqalignno{
[x^n]\ y &= {1\over n} [x^{n-1}]\ \phi(x)^n \cr
&= {1\over n} [x^{n-1}]\ {\prod_{m\geq 1} (1-y^m)^{-n}} \cr
&= {1\over n} [x^{n-1}]\ {\sum_{\l\in\setP} \prod_{v\in\l} 
\bigl( 1+{n-1\over h_v^2}\bigr)x
} \cr
&= {1\over n}
\sum_{\l\vdash n-1}\ \prod_{v\in\l}
\bigl(1+{n-1\over h_v^2} \bigr).\qed \cr
}
$$

As the coefficients of $y(x)$ are all positive integers we have the following
result.

\proclaim Corollary \sec.2 [=1.6].
For any positive integer $n$ the following expression
$$
{1\over n+1}\sum_{\l\vdash n} \prod_{v\in\l} \bigl(1+{n\over h_v^2}\bigr)
$$
is a positive integer.


\bigskip
\medskip
{\bf Acknowledgements}.
The author wishes to thank Dominique Foata for 
helpful discussions during the preparation of this paper.
He also thank Alain Lascoux, Richard Stanley and Kathy Ji for comments 
on a previous version.

\vskip 4mm 
\medskip

\bigskip \bigskip


\centerline{References}

{\eightpoint

\bigskip 
\bigskip

\divers AF02|Adin, Ron M.; Frumkin, Avital|
Rim Hook Tableaux and Kostant's $\eta$-Function Coefficients,
{\it arXiv: math.CO/0201003}| 

\livre An76|Andrews, George E.|The Theory of
Partitions|Addison-Wesley, Reading, {\oldstyle 1976}
({\sl Encyclopedia of Math. and Its Appl.,} vol.~{\bf 2})|


\article Be98|Bessenrodt, Christine|On hooks of Young diagrams|Ann. of 
Comb.|2|1998|103--110|

\divers BM02|Bacher, Roland; Manivel, Laurent|Hooks and Powers of Parts in 
Partitions, 
{\sl S\'em.  Lothar. Combin.}, vol.~{\bf 47}, article B47d, {\oldstyle 2001}, 
11 pages|

\divers BG06|Berkovich, Alexander; Garvan, Frank G.|%
The BG-rank of a partition and its applications, 
{\it arXiv: math/0602362}|

\divers CFP05|Cellini, Paola; Frajria, Pierluigi M.; Papi, Paolo|%
The $\hat W$-orbit of $\rho$, Kostant's formula for powers of the Euler 
product and affine Weyl groups as permutations of $\setZ$, 
{\it arXiv: math.RT/0507610}|

\divers CJW08|Chen, William Y.C.; Ji, Kathy Q.; Wilf, Herbert S.|BG-ranks 
and 2-cores, {\it arXiv: math/0605474v2}, {\oldstyle 2008}|

\article Dy72|Dyson, Freeman J.|Missed opportunities|%
Bull. Amer. Math. Soc.|78|1972|635--652|

\divers Eu83|Euler, Leonhard|The expansion of the infinite product 
$(1-x)(1-xx)(1-x^3)(1-x^4)(1-x^5)(1-x^6)$ etc. into a single series, 
{\sl English translation from the Latin by Jordan Bell} 
on {\it arXiv:math.HO/0411454}|

\article FK99|Farkas, Hershel M.; Kra, Irwin|On the Quintuple Product
Identity|Proc. Amer. Math. Soc.|27|1999|771--778|

\divers FH99|Foata, Dominique; Han, Guo-Niu|%
The triple, quintuple and septuple product identities revisited. 
{\sl Sem. Lothar. Combin.} Art. B42o, 12 pp|

\article FRT54|Frame, J. Sutherland; Robinson, Gilbert de Beauregard;        
Thrall, Robert M.|The hook graphs of the symmetric groups|Canadian 
J. Math.|6|1954|316--324|

\article GKS90|Garvan, Frank; Kim, Dongsu; Stanton, Dennis|Cranks and 
$t$-cores|Invent. Math.|101|1990|1--17|

\article GNW79|Greene, Curtis; Nijenhuis, Albert;
Wilf, Herbert S.|A probabilistic proof of a formula for the number of     
Young tableaux of a given shape|Adv. in Math.|31|1979|104--109|

\article Ho86|Hoare, A. Howard M.|An Involution of Blocks in the Partitions 
of $n$|Amer. Math. Monthly|93|1986|475--476|

\article JS89|Joichi, James T.; Stanton, Dennis|An
involution for Jacobi's identity|Discrete Math.|73|1989|261--271|

\article Ka74|Kac, Victor G.|Infinite-dimensional Lie algebras and 
Dedekind's $\eta$-function|Functional Anal. Appl.|8|1974|68--70|

\article Kn70|Knuth, Donald E.|Permutations, matrices, and generalized Young
tableaux|Pacific J. Math.|34|1970|709-727|

\livre Kn98|Knuth, Donald E.|The Art of Computer Programming,  {\bf vol.~3}, 
Sorting and Searching, 2nd ed.|Addison Wesley Longman,  {\oldstyle 1998}|

\article Ko76|Kostant, Bertram|On Macdonald's $\eta$-function formula, the
Laplacian and generalized exponents|Adv. in Math.|20|1976|179--212|

\article Ko04|Kostant, Bertram|Powers of the Euler product and commutative 
subalgebras of a complex simple Lie algebra|Invent. Math.|158|2004|181--226|

\article Kr99|Krattenthaler, Christian|\ Another involution
principle-free bijective proof of  
Stanley's hook-content formula|J.      
Combin. Theory Ser. A|88|1999|66--92| 

\article KS82|Kirdar, M. S.; Skyrme, Tony H. R.|On an Identity Related to 
Partitions and Repetitions of Parts|Canad. J. Math.|34|1982|194-195|

\livre La01|Lascoux, Alain|Symmetric Functions and Combinatorial Operators on 
Polynomials|CBMS Regional Conference Series in Mathematics, Number 99, 
{\oldstyle 2001}|

\article Le78|Lepowsky, James|Macdonald-type identities|%
Advances in Math.|27|1978|230--234| 

\article Ma72|Macdonald, Ian G.|\ Affine root systems and 
Dedekind's $\eta $-function|Invent. Math.|15|1972|91--143|

\livre Ma95|Macdonald, Ian G.|Symmetric Functions and Hall Polynomials|
Second Edition, Clarendon Press, Oxford, {\oldstyle 1995}|

\article Mi85|Milne, Stephen C.|An elementary proof of
the Macdonald identities for $A\sp {(1)}\sb l$|Adv. in Math.|
57|1985|34--70|

\article Mo75|Moody, Robert V.|\ Macdonald identities and Euclidean 
Lie algebras|Proc. Amer. Math. Soc.|48|1975|43--52|

\article NPS97|Novelli, Jean-Christophe; Pak, Igor;
Stoyanovskii, Alexander V.|A direct bijective proof of the hook-length
formula|Discrete Math. Theor. Comput. Sci.|1|1997|53--67|

\article Ro58|Robinson, G. de B.|A remark by Philip Hall|Canad. 
Math. Bull.|1|1958|21--23|

\article RS06|Rosengren, Hjalmar; Schlosser, Michael|%
Elliptic determinant evaluations and the Macdonald identities 
for affine root systems|Compositio Math.|142|2006|937-961|

\article RW83|Remmel, Jeffrey B.; Whitney, Roger A|
bijective proof of the hook formula for the number of column strict tableaux
with bounded entries|European J. Combin.|4|1983|45--63|

\livre Se70|Serre, Jean-Pierre|Cours d'arithm\'etique|Collection SUP: 
``Le Math\'emati\-cien", 2 Presses Universitaires de France, Paris 
{\oldstyle 1970}|

\divers Slo|Sloane, Neil; {\it al.}|
The On-Line Encyclopedia of Integer Sequences, {
\tt http:// www.research.att.com/\char126njas/sequences/}|

\livre St99|Stanley, Richard P.|Enumerative Combinatorics, vol. 2|
Cambridge university press, {\oldstyle 1999}|

\divers St04|Stanley, Richard P.|{\sl Errata and Addenda to Enumerative 
Combinatorics Volume 1, Second Printing}, Rev. Feb. 13, 2004.
{\tt http://www-math.mit.edu/\char126rstan/ ec/newerr.ps}|

\divers Ve|Verma, Daya-Nand|Review of the paper ``Affine root systems and 
Dedekind's $\eta $-function" written by Macdonald, I. G., MR0357528(50\#9996),
MathSciNet, 7 pages|

\divers We1|Weisstein, Eric W.|Elder's Theorem, from MathWorld -- 
A Wolfram Web Resource|

\divers We2|Weisstein, Eric W.|Stanley's Theorem, from MathWorld -- 
A Wolfram Web Resource|

\article Wi69|Winquist, Lasse|%
An elementary proof of $p(11m+6)\equiv 0\,({\rm mod} 11)$|%
J. Combinatorial Theory|6|1969|56--59|

\article Ze84|Zeilberger, Doron|A short hook-lengths bijection
inspired by the Greene-Nijen\-huis-Wilf proof|Discrete 
Math.|51|1984|101--108|

\bigskip

\irmaaddress
}
\vfill\eject

\end